\documentclass[a4paper,12pt]{article}
\usepackage{amsmath, amsthm}
\usepackage{amssymb,amsfonts}
\usepackage{enumerate}
\usepackage{csquotes}
\usepackage{fancyhdr}
\usepackage{graphicx}
\usepackage{hyperref}
\usepackage{authblk}

\usepackage{amssymb}

\theoremstyle{definition}
\newtheorem{definition}{Definition}[section]

\newtheorem{remark}[definition]{Remark}

\theoremstyle{plain}
\newtheorem{theorem}[definition]{Theorem}
\newtheorem{lemma}[definition]{Lemma}

\newtheorem{corollary}[definition]{Corollary}
\numberwithin{equation}{section}

\usepackage{color}   %May be necessary if you want to color links
\usepackage{hyperref}
\usepackage{url}
\begin{document}
\date{}
\author[1]{Dr. Rupali S. Jain}
\author[2]{Dr. B. Surendranath Reddy}
\author[3]{Mr. Wajid M. Shaikh}
\affil[1]{Associate Professor, School of Mathematical Sciences, S.R.T.M. University, Nanded}
\affil[2]{Assistant Professor, School of Mathematical Sciences, S.R.T.M. University, Nanded}
\affil[3]{Research Scholar, School of Mathematical Sciences, S.R.T.M. University, Nanded}
\title{
    {\bf Construction of Linear Codes from the Unit Graph $G(\mathbb{Z}_{n})$}\\
}
\maketitle
\begin{abstract}
In this paper, we  consider the unit graph $G(\mathbb{Z}_{n})$, where $n=p_{1}^{n_{1}} \text{ or } p_{1}^{n_{1}}p_{2}^{n_{2}} \text{ or } p_{1}^{n_{1}}p_{2}^{n_{2}}p_{3}^{n_{3}}$ and $p_{1}, p_{2}, p_{3}$ are distinct primes.  For any prime $q$, we construct $q$-ary linear codes from the incidence matrix of the unit graph $G(\mathbb{Z}_{n})$ with their parameters. We also prove that the dual of the constructed codes have minimum distance either 3 or 4. Lastly, we stated two conjectures on diameter of unit graph $G(\mathbb{Z}_{n})$ and  linear codes constructed from the incidence matrix of the unit graph $G(\mathbb{Z}_{n})$ for any integer $n$.
\end{abstract}
\section{Introduction}
In 1990, the unit graph was first introduced by Grimaldi, R. P.\cite{14} for $\displaystyle \mathbb{Z}_{n}$, is the graph with no loops and parallel edges and $\displaystyle x,y\in \mathbb{Z}_{n}$ are adjacent if and only if $x+y$ is a unit in $\mathbb{Z}_{n}$.
 \par In 2010, Fish, W., Key, J. D., $\&$ Mwambene, E.\cite{7}, constructed linear codes form incidence matrices of line graphs of Hamming graphs. Key, J. D., $\&$ Rodrigues, B. G. \cite{8} constructed  codes from lattice graphs and examine their decoding techniques using permutation decoding. Similar type of results were examined by several researcher \cite{9,10}. In 2013, Dankelmann, P., Key, J. D., $\&$ Rodrigues, B. G. \cite{3} gave the generalization of relationship between parameters of connected graphs and codes generated from incidence matrices of graphs and obtained upper bounds for minimum distance of their dual codes. Also several researchers carried out their research in construction of linear codes from adjacency matrices of  some special graphs, as mentioned in \cite{11,12,13}.
\par Recently, in 2021,  Annamalai, N., $\&$ Durairajan, C.\cite{1}, constructed linear codes from the incidence matrices of unit graph $\displaystyle G(\mathbb{Z}_{p})$ and $\displaystyle G(\mathbb{Z}_{2p})$, for odd prime $p$. In this paper, we generalise their work by constructing linear codes from the  incidence matrices of unit graph $\displaystyle G(\mathbb{Z}_{n})$, where $n=p_{1}^{n_{1}} \text{ or } p_{1}^{n_{1}}p_{2}^{n_{2}} \text{ or } p_{1}^{n_{1}}p_{2}^{n_{2}}p_{3}^{n_{3}}$ and $p_{1}, p_{2}, p_{3}$ are distinct primes. We also find the parameters  for these codes and their dual codes. Finally we conclude by stating two conjectures.
%\par In this article, section 2 is about preliminaries.
\section{Preliminaries}
In this section, we recall definitions and results related to  unit graphs and linear codes. Let $\mathbb{Z}_{n}$ denote the ring of integers modulo $n$. Here, we denote units and non-units of $\displaystyle \mathbb{Z}_{n}$ by $\displaystyle U(\mathbb{Z}_{n})$ and $\displaystyle N_{U}(\mathbb{Z}_{n})$ respectively.
\begin{definition}\cite{18}[Linear Code]
Let $\displaystyle \mathbb{F}_{q}$ represents the finite field with $q$ elements. A  linear code $C_{q}$ of length $n$ is a subspace of $\displaystyle \mathbb{F}^{n}_{q}$ and it is called $q$-ary linear code. Dimension of linear code $C_{q}$ is the dimension of $C_{q}$ as a vector space over field  $\displaystyle \mathbb{F}_{q}$ and is denoted by $\displaystyle \text{dim}(C_{q})$.
\end{definition}
\begin{definition}\cite{18}[Dual of code]
Let $C_{q}$ be a linear code of length $n$ over $\displaystyle \mathbb{F}_{q}$. Then dual of code $C_{q}$ is the orthogonal compliment of the subspace $C_{q}$ in $\displaystyle \mathbb{F}^{n}_{q}$ and is denoted by $\displaystyle C^{\perp}_{q}$.
\end{definition}
\begin{theorem}\cite{18}\label{thm8}
Let $C_{q}$ be a $q$-ary code of length $n$ over a field $\displaystyle \mathbb{F}_{q}$. Then $\displaystyle C^{\perp}_{q}$ is a linear code of length $n$ and $\displaystyle \text{dim}(C_{q}^{\perp})=n-\text{dim}(C_{q})$.
\end{theorem}
\begin{definition}
\begin{enumerate}
  \item Let $x$ and $y$ be vectors in $\displaystyle \mathbb{F}^{n}_{q}$. Then Hamming distance of $x$ and $y$, denoted by $\displaystyle d_{C}(x,y)$, is defined to be the number of places at which $x$ and $y$ differ.
  \item Let $x$ be a vector in $\displaystyle \mathbb{F}^{n}_{q}$. Then Hamming weight of $x$ is defined to be the number of non-zero coordinates in $x$ and is denoted by $\text{wt}(x)$. Clearly, $\displaystyle \text{wt}(x-y)=d_{C}(x,y)$.
\end{enumerate}
\end{definition}

\begin{definition}\cite{18}[Minimum Hamming weight]
Let $C_{q}$ be a linear code. Then minimum Hamming weight of $C_{q}$, denoted by $\text{wt}(C_{q})$, is defined as
$\displaystyle \text{wt}(C_{q})=\text{min}\{\text{wt}(x) \ | \ x\in C_{q} \ \& \ x\neq 0\}$.
\end{definition}

\begin{remark}\cite{18}[Minimum Hamming distance]
Let $C_{q}$ be a linear code. The minimum Hamming distance of code $C_{q}$, denoted by $d(C_{q})$, is defined as
$\displaystyle d(C_{q})=\text{min}\{d_{C}(x,y) \ | \ x,y\in C_{q} \ \& \ x\neq y\}$. Note that $d(C_{q})=\text{wt}(C_{q})$.
\end{remark}

\begin{remark}\cite{18}
A $q$-ary linear code $C_{q}$ of length $n$, dimension $k$ and minimum distance $d$ is called $[n,k,d]_{q}$ linear code.
\end{remark}
\begin{definition}\cite{18}
A generator matrix of linear code $C_{q}$ is a matrix $H$ whose rows form a basis for $C_{q}$ and a generator matrix $\displaystyle H^{\perp}$ of linear code $C^{\perp}_{q}$ is called parity-check matrix of $C_{q}$.
\end{definition}
\par Let $\displaystyle G=(V,E)$ be a graph with vertex set $V$  and edge set $E$. For any $x,y\in V$, $[x,y]$ denote the edge between $x$ and $y$ and if $[x,y]\in E$ then we call $x$ is adjacent to $y$. A graph is called simple if it does not have loops and parallel edges. A complete graph is a simple graph in which any distinct pair of vertices is joined by an edge. If the vertex set $V$ of $G$ can be partitioned into two non-empty subsets $W_{1}$ and $W_{2}$ such that  each edge in $G$ has one end in $W_{1}$ and one end in $W_{2}$ then $G$ is called bipartite.
A complete bipartite graph is a simple bipartite graph $G$, with bipartition $V=W_{1}\cup W_{2}$, in which every vertex in $W_{1}$ is joined to every vertex in $W_{2}$.
\par An edge $e$ of graph $G$ is said to be incident with the vertex $x$, if $x$ is the end vertex of $e$. In this case we also say that $x$ is incident with $e$. The degree of vertex $x$ is the number of edges of $G$ incident with $x$ and it is denoted by $\text{deg}(x)$. If $x\in V$ and $\text{deg}(x)\leq \text{deg}(y)$ for all $y\in V$ then $\text{deg}(x)$ is called minimum degree of $G$ and it is denoted by $\delta(G)$.
If for some positive integer $k$, $\text{deg}(x)=k$ for every vertex of the graph of $G$, then $G$ is called $k$-regular graph. A vertex $x$ is said to be connected to a vertex $y$ in a graph $G$ if there is path in $G$ from $x$ to $y$. A graph is called connected if every two vertices are connected. A nontrivial closed trail in graph $G$ is called cycle if its origin and internal vertices are distinct. A cycle of length $k$, i.e. with $k$ edges is called $k$-cycle. The girth of graph $G$, is denoted by $g_{r}(G)$, is the length of the shortest cycle contained in $G$.
\begin{definition}
The distance between two vertices $x$ and $y$, denoted by $d(x,y)$, is the length of a shortest path from $x$ to $y$. The diameter of a graph $G$ is denoted by $\text{diam}(G)$, is the maximum distance between any two vertices in $G$. i.e.
$\displaystyle \text{diam}(G)=\text{Max}\{d(x,y) \ | \ x,y\in V\}$.
\end{definition}

\begin{definition}\cite{19}
Let $G$ be a simple graph. The edge connectivity of $G$, denoted by $\lambda(G)$, is the smallest number of edges in $G$ whose deletion from $G$ either leaves a disconnected graph or an empty graph.
\end{definition}
\begin{definition}\cite{2}
	Let $R$ be a ring with nonzero identity. The unit graph of $R$, denoted by $G(R)$, is a graph with vertex set as $R$ and two distinct vertices $x$ and $y$ are adjacent if and only if $x + y$ is a unit of $R$.
\end{definition}
\begin{theorem}\label{thm3}\cite{4}
Let $\displaystyle G=(V,E)$ be a connected graph with vertex set $V$.\\ If $\text{diam}(G)\leq 2$ then edge connectivity of $G$ is  $\lambda(G)=\delta(G)$.
\end{theorem}

\begin{theorem}\label{thm4}\cite{5}
Let $\displaystyle G=(V,E)$ be a connected bipartite graph, if $\displaystyle \text{diam}(G)\leq 3$ then edge connectivity of $G$ is $\displaystyle \lambda(G)=\delta(G)$.
\end{theorem}

\begin{theorem}\label{thm1}\cite{2}
Let $\mathcal{R}$ be a finite ring. Then the following statements hold for the unit graph of $\mathcal{R}$
\begin{enumerate}
  \item If $\displaystyle 2\notin U(\mathcal{R})$, then the unit graph $G(\mathcal{R})$ is a $|U(\mathcal{R})|$-regular graph.
  \item If $\displaystyle 2\in U(\mathcal{R})$, then for every $x\in U(\mathcal{R})$ we have $\displaystyle \text{deg}(x)=|U(\mathcal{R})|-1$ and for every $\displaystyle x\in N_{U}(\mathcal{R})$ we have $\displaystyle \text{deg}(x)=|U(\mathcal{R})|$.
\end{enumerate}
\end{theorem}

\begin{theorem}\label{thm7}\cite{15}
Let $R$ be a ring. Then
\begin{enumerate}
  \item If $\displaystyle |U(R)|= 2$ then $g_{r}(G)\in \{3,4,6\}$.
  \item If $\displaystyle |U(R)|\geq 3$ then $g_{r}(G)\in \{3,4\}$.
\end{enumerate}
\end{theorem}

\begin{theorem}\label{thm2}\cite{3}
Let $G=(V,E)$ be a connected graph and let $H$ be a $\displaystyle |V|\times |E|$ incidence matrix for $G$. Then binary code generated by $H$ is\\ $\displaystyle C_{2}(H)=\displaystyle [|E|,|V|-1,\lambda(G)]_{2}$.
\end{theorem}

\begin{theorem}\cite{3}\label{thm5}
  Let $G=(V,E)$ be a connected bipartite graph and let $H$ be a $\displaystyle |V|\times |E|$ incidence matrix for $G$, and $q$ be an odd prime. Then $q$-ary code generated by $H$ is $\displaystyle C_{q}=[|E|,|V|-1,\lambda(G)]_{q}$.
\end{theorem}

\begin{theorem}\cite{3}\label{thm6}
Let $G$ be a connected graph with girth $g_{r}(G)$ and even girth $g_{r}(G)_{e}$. Let $H$ be an incidence matrix for $G$, $C=C_{q}(H)$, where $q$ is any prime, and $\displaystyle d^{\perp}$ is the minimum distance of $\displaystyle C_{q}^{\perp}$. If $q=2$ or $g_{r}(G)$ is even then $\displaystyle d^{\perp}=g_{r}(G)$.
\end{theorem}

\section{Construction of Linear Codes from the Incidence Matrix of a Unit Graph $G(\mathbb{Z}_{p^{n}})$}
In this section, we construct binary and $q$-ary linear codes $C_{2}$ and $C_{q}$ generated from the incidence matrix of the unit graph $\displaystyle G(\mathbb{Z}_{p^{n}})$, where $p$ is any prime number and $n\in\mathbb{N}$. We also examine the dual codes $C^{\perp}_{2}$ and $C^{\perp}_{q}$ with their parameters.
\par
Let $\displaystyle \mathbb{Z}_{p^{n}}$ denote the ring. Then $\displaystyle U(\mathbb{Z}_{p^{n}})=\{x\in \mathbb{Z}_{p^{n}} \ | \ \text{g.c.d}(x,p^{n})=1\}$ and $\displaystyle N_{U}(\mathbb{Z}_{p^{n}})=\{x\in \mathbb{Z}_{p^{n}} \ | \ x=\alpha p \}\neq \phi$. Let $G(\mathbb{Z}_{p^{n}})$ be a unit graph with vertex set $\displaystyle V=\mathbb{Z}_{p^{n}}$ and $x,y\in \mathbb{Z}_{p^{n}}$ are adjacent if and only if $x+y\in U(\mathbb{Z}_{p^{n}})$.
\begin{theorem}\label{thm21}
Let $\displaystyle G(\mathbb{Z}_{p^{n}})$ be a unit graph and $p$ is an odd prime. Then the graph $\displaystyle G(\mathbb{Z}_{p^{n}})$ is connected with $\displaystyle |V|=p^{n}$
and $\displaystyle |E|=\frac{(p^{n}-1)\phi(p^{n})}{2}$.
\begin{proof}
Clearly, $\displaystyle |V|=p^{n}$. Note that for every $\displaystyle x\in U(\mathbb{Z}_{p^{n}})$ and $\displaystyle y\in N_{U}(\mathbb{Z}_{p^{n}})$, $[x,y]$ is an edge in $G(\mathbb{Z}_{p^{n}})$. Suppose this is not true which gives $x+y$ is not a unit, this implies $\displaystyle p|x+y$. But $\displaystyle y\in N_{U}(\mathbb{Z}_{p^{n}})$ i.e. $y=\alpha p$ this implies $p|x $ which is a contradiction. Hence $\displaystyle G(\mathbb{Z}_{p^{n}})$ is connected graph.
\par Since $p$ is odd prime, $2\in U(\mathbb{Z}_{p^{n}})$. By theorem (\ref{thm1}), we have $\displaystyle\text{deg}(x)=|U(\mathbb{Z}_{p^{n}})|-1=\phi(p^{n})-1, \ \forall \ x\in U(\mathbb{Z}_{p^{n}})$ and $\displaystyle \text{deg}(x)=\phi(p^{n}), \ \forall \ x\in N_{U}(\mathbb{Z}_{p^{n}})$. Now
\begin{align*}
   |E| & =\frac{\sum_{x\in \mathbb{Z}_{p^{n}}}\text{deg}(x)}{2} \\
       & =\frac{\sum_{x\in U(\mathbb{Z}_{p^{n}})}\text{deg}(x)+\sum_{x\in N_{U}(\mathbb{Z}_{p^{n}})}\text{deg}(x)}{2} \\
                  & = \frac{\phi(p^{n})[\phi(p^{n})-1]+[p^{n}-\phi(p^{n})]\phi(p^{n})}{2}\\
                  & = \frac{\phi(p^{n})[\phi(p^{n})-1+p^{n}-\phi(p^{n})]}{2}\\
   |E| & = \frac{\phi(p^{n})[p^{n}-1]}{2}
\end{align*}
\end{proof}
\end{theorem}

\begin{theorem}\label{thm22}
 Let $\displaystyle G(\mathbb{Z}_{p^{n}})$ be a unit graph and $p$ is an odd prime. Then edge connectivity of $\displaystyle G(\mathbb{Z}_{n})$ is $\displaystyle \lambda(G(\mathbb{Z}_{p^{n}}))=\displaystyle \phi(p^{n})-1$.
\begin{proof}
First we show that $\text{diam}(G(\mathbb{Z}_{p^{n}}))\leq 2$. For any $\displaystyle x,y\in \mathbb{Z}_{p^{n}}$ we have following cases,\\
{\bf Case I:} Let $\displaystyle x,y\in U(\mathbb{Z}_{p^{n}})$. Then there exist $\displaystyle z\in N_{U}(\mathbb{Z}_{p^{n}})$ such that $x$ and $y$ adjacent to $z$. Hence $d(x,y)\leq 2$.\\
{\bf Case II:} If $\displaystyle x,y\in N_{U}(\mathbb{Z}_{p^{n}})$ then for any $z\in U(\mathbb{Z}_{p^{n}})$, we have $[x,z]$ and $[z,y]$ are edges in $\displaystyle G(\mathbb{Z}_{p^{n}})$, which gives $d(x,y)\leq 2$\\
{\bf Case III:} If $\displaystyle x\in U(\mathbb{Z}_{p^{n}})$ and $\displaystyle y\in N_{U}(\mathbb{Z}_{p^{n}})$ then $x$ is adjacent to $y$, which gives $d(x,y)= 1$.\\
We get $d(x,y)\leq 2$ for all $\displaystyle x,y\in \mathbb{Z}_{p^{n}}$ which gives $\text{diam}(G(\mathbb{Z}_{p^{n}}))\leq 2$.\\
 Now by Theorem \ref{thm3}, we have $\displaystyle \lambda(G(\mathbb{Z}_{p^{n}}))=\delta(G(\mathbb{Z}_{p^{n}}))=\phi(p^{n})-1$.
\end{proof}
\end{theorem}

\begin{corollary}\label{thm26}
Let $\displaystyle G(\mathbb{Z}_{p^{n}})$ be a unit graph, where $p$ is an odd prime and $\displaystyle p^{n}\neq 3$. Then $\displaystyle g_{r}(G(\mathbb{Z}_{p^{n}}))=3$.
\begin{proof}
If $n=1$ then $p\geq3$. Note that $\displaystyle 1,2\in \mathbb{Z}_{p}$ and $1$ is adjacent to $2$. Since, $\displaystyle 0\in \mathbb{Z}_{p}$ which gives $1$ and $2$ are adjacent to $0$. Hence, we get cycle of length 3, which conclude our result. If $n>1$ then for any $\displaystyle x\in U(\mathbb{Z}_{p^{n}})$ we have $\displaystyle x+p\in U(\mathbb{Z}_{p^{n}})$. Also $x$ is adjacent to $x+p$. Note that  $x$ and $x+p$ are adjacent to $p$. Hence, we get cycle of length $3$, which implies $g_{r}(G(\mathbb{Z}_{p^{n}}))=3$.
\end{proof}
\end{corollary}

\begin{lemma}\label{thm23}
Let $\displaystyle G(\mathbb{Z}_{2^{n}})$ be a unit graph. Then $x$ and $y$ are adjacent if and only if $x\in U(\mathbb{Z}_{2^{n}})$ and $y\in N_{U}(\mathbb{Z}_{2^{n}})$.
\begin{proof}
Let $G(\mathbb{Z}_{2^{n}})$ be a unit graph.
Suppose that, $x$ and $y$ are adjacent for some $x,y\in \mathbb{Z}_{2^{n}}$.  That is $[x,y]$ is an edge in $G(\mathbb{Z}_{2^{n}})$, which implies $x+y\in U(\mathbb{Z}_{2^{n}})$, hence $x+y=2k+1$ for some integer $k$. Assume that $\displaystyle x,y\in U(\mathbb{Z}_{2^{n}})$ which imlies $x=2k_{1}+1$ and $y=2k_{2}+1$ for some integer $k_{1}$ and $k_{2}$. From this, we have  $\displaystyle 2|x+y$,  which gives a contradiction to our assumption. Similarly, if $\displaystyle x,y\in N_{U}(\mathbb{Z}_{p^{n}})$, then $\displaystyle 2|x+y$, which is again a contradiction. Hence $\displaystyle x\in U(\mathbb{Z}_{p^{n}}) $ and $\displaystyle y\in N_{U}(\mathbb{Z}_{2^{n}})$.\\
Conversely, let $\displaystyle x\in U(\mathbb{Z}_{2^{n}})$  and $\displaystyle y\in N_{U}(\mathbb{Z}_{2^{n}})$. Then $\displaystyle x=2k_{1}$ and $\displaystyle y=2k_{1}+1$ this implies $2\nmid x+y$ and hence $[x,y]$ is an edge in $G(\mathbb{Z}_{2^{n}})$.
\end{proof}
\end{lemma}

\begin{corollary}\label{thm24}
Let $\displaystyle G(\mathbb{Z}_{2^{n}})$ be a unit graph. Then $\displaystyle G(\mathbb{Z}_{2^{n}})$ is complete bipartite graph with bipartition $\displaystyle W_{1}=\{x\in \mathbb{Z}_{2^{n}} \ | \ x=2\alpha\ \text{for some }\ \alpha\in\mathbb{Z} \}$ and $\displaystyle W_{2}=\{x\in \mathbb{Z}_{2^{n}} \ | \ x=2\alpha+1\ \text{for some }\ \alpha\in\mathbb{Z} \}$.
\end{corollary}

\begin{theorem}\label{thm25}
Let $\displaystyle G(\mathbb{Z}_{2^{n}})$ be a unit graph. Then
\begin{enumerate}
  \item $\displaystyle |V|=2^{n}$ and $\displaystyle |E|=2^{2(n-1)}$.
  \item $\displaystyle \lambda(G(\mathbb{Z}_{2^{n}}))=\phi(2^{n})$.
\end{enumerate}
\begin{proof}
\begin{enumerate}
  \item  By the definition of unit graph $\displaystyle G(\mathbb{Z}_{2^{n}})$, we have $\displaystyle |V|=2^{n}$. By Corollary (\ref{thm24}), $\displaystyle G(\mathbb{Z}_{2^{n}})$ is complete bipartite graph, which implies $\displaystyle |E|=2^{n-1}2^{n-1}=2^{2(n-1)}$.
  \item From Corollary (\ref{thm24}), $\displaystyle G(\mathbb{Z}_{2^{n}})$ is complete bipartite graph and hence $\displaystyle \lambda(G(\mathbb{Z}_{2^{n}}))=\phi(2^{n})$.
\end{enumerate}
\end{proof}
\end{theorem}
\begin{theorem}
Let $\displaystyle G(\mathbb{Z}_{p^{n}})$ be a unit graph and $H$ be a $|V|\times |E|$ incidence matrix of $G(\mathbb{Z}_{p^{n}})$.
\begin{enumerate}
  \item If $p$ is an odd prime then binary code generated by $H$ is a \\ $\displaystyle C_{2}(H)=\left[\frac{\phi(p^{n})(p^{n}-1)}{2},p^{n}-1,\phi(p^{n})-1\right]_{2}$ code over the finite field $\displaystyle \mathbb{F}_{2}$.
  \item If $p=2$ then for any odd prime $q$, the $q$-ary code generated by $H$ is a $\displaystyle C_{q}(H)=[2^{2(n-1)},2^{n}-1,2^{(n-1)}]_{q}$ code over the finite field $\displaystyle \mathbb{F}_{q}$.
\end{enumerate}
\begin{proof}
\begin{enumerate}
  \item Let $G(\mathbb{Z}_{p^{n}})$ be a unit graph, where $p$ is an odd prime and $H$ be an incidence matrix of $G(\mathbb{Z}_{p^{n}})$.  By Theorem (\ref{thm21}),  $G(\mathbb{Z}_{p^{n}})$ is connected graph and hence by Theorem (\ref{thm2}), binary code generated by $H$ is $C_{2}(H)=[|E|,|V|-1,\lambda(G(\mathbb{Z}_{p^{n}}))]_{2}$. By theorem (\ref{thm21}) and (\ref{thm22}), we get
$\displaystyle |E|=\frac{\phi(p^{n})(p^{n}-1)}{2}, \ |V|=p^{n}$  and the  edge connectivity of $\displaystyle G(\mathbb{Z}_{p^{n}})$ is $\lambda(G(\mathbb{Z}_{p^{n}}))=\phi(p^{n})-1$. Hence we get $\displaystyle C_{2}(H)=\left[\frac{\phi(p^{n})(p^{n}-1)}{2}, p^{n}-1, \phi(p^{n})-1\right]_{2}$.
\item Let $p=2$. Then  by Corollary (\ref{thm24}), $\displaystyle G(\mathbb{Z}_{2^{n}})$ is a complete bipartite graph, which implies $G(\mathbb{Z}_{2^{n}})$ is connected bipartite graph and hence by Theorem (\ref{thm2}) and (\ref{thm5}), for any prime $q$, $q$-ary code generated by $H$ is $C_{q}(H)=[|E|,|V|-1,\lambda(G(\mathbb{Z}_{2^{n}}))]_{q}$. Using Theorem (\ref{thm25}), we get
    $\displaystyle |E|=2^{2(n-1)}, \ |V|=2^{n}$  and the  edge connectivity of $\displaystyle G(\mathbb{Z}_{2^{n}})$ is $\lambda(G(\mathbb{Z}_{2^{n}}))=2^{n-1}$. Hence $\displaystyle C_{q}(H)=[2^{2(n-1)}, 2^{n}-1, 2^{n-1}]_{q}$.
\end{enumerate}
\end{proof}
\end{theorem}

\begin{corollary}
Let $C_{q}(H)$ and $C_{2}(H)$ denote the codes generated by incidence matrix of $\displaystyle G(\mathbb{Z}_{2^{n}})$ and $\displaystyle G(\mathbb{Z}_{p^{n}})$. Then
\begin{enumerate}
  \item Dual of code $C_{2}(H)$ is $\displaystyle C^{\perp}_{2}=\left[ \frac{(p^{n}-1)\phi(p^{n})}{2}, \frac{(p^{n}-1)[\phi(p^{n})-2]}{2}, 3
   \right]_{2}$, where $p^{n}\neq 3$.
   \item Dual of code $C_{q}(H)$ is $\displaystyle C^{\perp}_{q}=\left[2^{2(n-1)},2^{n}(2^{n-2}-1)+1,4\right]_{q}$.
\end{enumerate}
\begin{proof}
\begin{enumerate}
  \item It follows from Theorems (\ref{thm8}) \&  (\ref{thm6}) and Corollary (\ref{thm26}).
  \item Since $\displaystyle G(\mathbb{Z}_{2^{n}})$ is complete bipartite graph, we have $\displaystyle g_{r}(G)=4$ and hence result follows from Theorems (\ref{thm8}) and (\ref{thm6}).
\end{enumerate}
\end{proof}
\end{corollary}

\section{Construction of Linear Codes from the Incidence Matrix of a Unit Graph $G(\mathbb{Z}_{p^{n_{1}}_{1}p^{n_{2}}_{2}})$}
In this section, we construct binary and $q$-ary linear codes from the incidence matrix of the unit graphs $\displaystyle G(\mathbb{Z}_{p^{n_{1}}_{1}p^{n_{2}}_{2}})$, where $p_{1}$ and $p_{2}$  are distinct primes and $\displaystyle n_{1},n_{2}\in \mathbb{N}$. We also discuss the parameters their dual codes.
\begin{lemma}\label{thm31}
	Let $\displaystyle G(\mathbb{Z}_{n})$ be a unit graph. If $2\in N_{U}(\mathbb{Z}_{n})$, then $\displaystyle G(\mathbb{Z}_{n})$ is bipartite.
	\begin{proof}
		We have  $\displaystyle W_{1}\cup W_{2}=\mathbb{Z}_{n}$, where $\displaystyle W_{1}=\{x\in \mathbb{Z}_{n} \ | \ x=2\alpha \ \text{for some } \alpha\in \mathbb{Z}\}$ and $\displaystyle W_{2}=\{x\in \mathbb{Z}_{n} \ | \ x=2\alpha+1 \ \text{for some } \alpha\in \mathbb{Z}\}$.  Now for any $\displaystyle x,y\in W_{1}$, we have $\displaystyle \text{gcd}(x+y,n)\geq 2$, which implies that $x$ is not adjacent to $y$. Similarly, for any $\displaystyle x,y\in W_{2}$, we have $\displaystyle \text{gcd}(x+y,n)\geq 2$. Hence $\displaystyle G(\mathbb{Z}_{n})=G(W_{1}\cup W_{2})$ is a bipartite graph.
	\end{proof}
\end{lemma}
\begin{theorem}\label{thm32}
Let $\displaystyle G(\mathbb{Z}_{p^{n_{1}}_{1}p^{n_{2}}_{2}})$ be a unit graph, where $\displaystyle p_{1}$ and  $\displaystyle p_{2}$ are distinct odd primes. Then $G(\mathbb{Z}_{p^{n_{1}}_{1}p^{n_{2}}_{2}})$ is a connected graph and $\displaystyle \text{diam}(G(\mathbb{Z}_{p^{n_{1}}_{1}p^{n_{2}}_{2}}))\leq 2$.
\begin{proof}
 We can rewrite
 \begin{align*}
 \mathbb{Z}_{p^{n_{1}}_{1}p^{n_{2}}_{2}}&=U(\mathbb{Z}_{p^{n_{1}}_{1}p^{n_{2}}_{2}})\cup N_{U}(\mathbb{Z}_{p^{n_{1}}_{1}p^{n_{2}}_{2}})\\
 & = U(\mathbb{Z}_{p^{n_{1}}_{1}p^{n_{2}}_{2}})\cup N_{p_{1}}\cup N_{p_{2}}\cup N_{p_{1}p_{2}}
 \end{align*}
 where
 \begin{align*}
 N_{p_{1}}&=\{x\in \mathbb{Z}_{p^{n_{1}}_{1}p^{n_{2}}_{2}} \ | \ x=\alpha p_{1} \ \& \ p_{2}\nmid \alpha \}\\
 N_{p_{2}}&=\{x\in \mathbb{Z}_{p^{n_{1}}_{1}p^{n_{2}}_{2}} \ | \ x=\alpha p_{2} \ \& \ p_{1}\nmid \alpha \}\\
 N_{p_{1}p_{2}}&=\{x\in \mathbb{Z}_{p^{n_{1}}_{1}p^{n_{2}}_{2}} \ | \ x=\alpha p_{1}p_{2}  \}
\end{align*}
 For $\displaystyle x,y\in \mathbb{Z}_{p^{n_{1}}_{1}p^{n_{2}}_{2}}$, we have  following cases\\
{\bf Case I:} If $\displaystyle x,y\in U(\mathbb{Z}_{p^{n_{1}}_{1}p^{n_{2}}_{2}})$, then $\displaystyle [x,0]$ and $\displaystyle [0,y]$ are edges in $\displaystyle G(\mathbb{Z}_{p^{n_{1}}_{1}p^{n_{2}}_{2}})$, which gives $d(x,y)\leq 2$.\\
{\bf Case II:} If $\displaystyle x\in U(\mathbb{Z}_{p^{n_{1}}_{1}p^{n_{2}}_{2}})$ and $\displaystyle y\in N_{U}(\mathbb{Z}_{p^{n_{1}}_{1}p^{n_{2}}_{2}})$, then we have the following possibilities  for $y$:\\
(a) If $\displaystyle y\in N_{p_{1}}$ ,then $y=\alpha p_{1}$. Suppose that $\displaystyle [x,y]=[x,\alpha p_{1}]$ is not an edge in $\displaystyle G(\mathbb{Z}_{p^{n_{1}}_{1}p^{n_{2}}_{2}})$ which implies $\displaystyle p_{2}\mid x+\alpha p_{1}$ i.e. $x+\alpha p_{1}=\beta p_{2}$. Note that $p_{1}\nmid \beta$, hence  $\displaystyle \beta p_{2}\in N_{p_{2}}$. Now we claim that $\displaystyle [x,\beta p_{2}]$ is an edge in $\displaystyle G(\mathbb{Z}_{p^{n_{1}}_{1}p^{n_{2}}_{2}})$. Suppose this is no true. Then  $p_{1}\mid x+\beta p_{2}$,  we get $\displaystyle p_{1}\mid 2x+\alpha p_{1}$, i.e. $\displaystyle p_{1}\mid x$ which is a contradiction. Hence, we get $x$ is adjacent to $\beta p_{2}$. Note that $\displaystyle [\beta p_{2},\alpha p_{1}]$ is an edge in $\displaystyle G(\mathbb{Z}_{p^{n_{1}}_{1}p^{n_{2}}_{2}})$. i.e. $\displaystyle [x,\beta p_{2}]$ and $\displaystyle [\beta p_{2},\alpha p_{1}]$ are edges in $\displaystyle G(\mathbb{Z}_{p^{n_{1}}_{1}p^{n_{2}}_{2}})$.  Hence $d(x,y)\leq 2$.\\
(b) If $\displaystyle y\in N_{p_{2}}$, then  $y=\alpha p_{2}$.  By following the above procedure, we get $d(x,y)\leq 2$.\\
(c) If $\displaystyle y\in N_{p_{1}p_{2}}$, then $y=\alpha p_{1}p_{2}$. Since $x$ is a unit, we get $\displaystyle [x,\alpha p_{1}p_{2}]$ is an edge in $\displaystyle G(\mathbb{Z}_{p^{n_{1}}_{1}p^{n_{2}}_{2}})$. For if this not hold, then we have $\displaystyle p_{1}\mid x+\alpha p_{1}p_{2}$ i.e. $p_{1}\mid x$ which is not possible. Similarly if $\displaystyle p_{2}\mid x+\alpha p_{1}p_{2}$ then $p_{2}\mid x$. Hence $d(x,y)=1\leq 2$.\\
{\bf Case III:} If $\displaystyle x,y\in N_{U}(\mathbb{Z}_{p^{n_{1}}_{1}p^{n_{2}}_{2}})$ then, consider the following possibilities\\
(a) If $\displaystyle x,y\in N_{p_{1}}$ then $x=\alpha p_{1}$ and $\displaystyle y=\alpha' p_{1}$. Since, $\displaystyle p_{2}\in N_{p_{2}}$ which implies $N_{p_{2}}\neq \phi$. Then for any $\displaystyle \beta p_{2}\in N_{p_{2}}$,. $\displaystyle [\alpha p_{1},\beta p_{2}]$ is an edges in $\displaystyle G(\mathbb{Z}_{p^{n_{1}}_{1}p^{n_{2}}_{2}})$. Similarly, we can prove that $\displaystyle [\beta p_{2},\alpha' p_{1}]$ is an edge in $\displaystyle G(\mathbb{Z}_{p^{n_{1}}_{1}p^{n_{2}}_{2}})$. i.e. $\displaystyle [x,\beta p_{2}]$ and $\displaystyle [\beta p_{2},y]$ are edges in $\displaystyle G(\mathbb{Z}_{p^{n_{1}}_{1}p^{n_{2}}_{2}})$. Hence, we have  $\displaystyle d(x,y)\leq 2$.\\
(b) If $\displaystyle x,y\in N_{p_{2}}$, then  $d(x,y)\leq2$ as in similar to case (a).\\
(c) If $\displaystyle x,y\in N_{p_{1}p_{2}}$, then $x=\alpha p_{1}p_{2}$ and $\displaystyle y=\alpha' p_{1}p_{2}$. Note that $\displaystyle U(\mathbb{Z}_{p^{n_{1}}_{1}p^{n_{2}}_{2}})\neq \phi$. Then for any $\displaystyle z\in U(\mathbb{Z}_{p^{n_{1}}_{1}p^{n_{2}}_{2}})$, we have  $\displaystyle [\alpha p_{1}p_{2},z]$ and $\displaystyle [z,\alpha' p_{1}p_{2}]$ are edges in $\displaystyle G(\mathbb{Z}_{p^{n_{1}}_{1}p^{n_{2}}_{2}})$. Hence $\displaystyle d(x,y)\leq 2$.\\
(d)If $\displaystyle x\in N_{p_{1}}$ and $\displaystyle y\in N_{p_{2}}$, then $x$ is adjacent to $y$, hence $\displaystyle d(x,y)=1\leq 2$.\\
(e) If $\displaystyle x\in N_{p_{1}}$ and $\displaystyle y\in N_{p_{1}p_{2}}$, then $x=\alpha p_{1} \ \& \ y=\alpha' p_{1}p_{2}$. Now, let $\beta p_{2}\in N_{2}$ then $\displaystyle z=\alpha p_{1}+\beta p_{2}\in U(\mathbb{Z}_{p^{n_{1}}_{1}p^{n_{2}}_{2}})$, from this we have $[z,\alpha' p_{1}p_{2}]$ is an edge. Also $\displaystyle [\alpha p_{1},z]$ is an edge in $\displaystyle G(\mathbb{Z}_{p^{n_{1}}_{1}p^{n_{2}}_{2}})$. Hence $d(x,y)\leq 2$.\\
(f) If $\displaystyle x\in N_{p_{2}}$ and $\displaystyle y\in N_{p_{1}p_{2}}$, then using the above procedure, we can show that $d(x,y)\leq2$.\\
Hence, $\displaystyle G(\mathbb{Z}_{p^{n_{1}}_{1}p^{n_{2}}_{2}})$  is connected graph and $d(x,y)\leq 2$, for all $\displaystyle x,y\in \mathbb{Z}_{p^{n_{1}}_{1}p^{n_{2}}_{2}}$, which gives $\text{diam}(G(\mathbb{Z}_{p^{n_{1}}_{1}p^{n_{2}}_{2}}))\leq 2$.
\end{proof}
\end{theorem}

\begin{corollary}\label{thm36}
Let $\displaystyle G(\mathbb{Z}_{p^{n_{1}}_{1}p^{n_{2}}_{2}})$ be a unit graph, where $\displaystyle p_{1}$ and  $\displaystyle p_{2}$ are odd primes. Then $\displaystyle g_{r}(G(\mathbb{Z}_{p^{n_{1}}_{1}p^{n_{2}}_{2}}))=3$.
\begin{proof}
Note that, $p_{1}$ is adjacent to $p_{2}$ and $p_{1}+p_{2}$ is adjacent to  both $p_{1} \ \& \ p_{2}$, from this we conclude the result.
\end{proof}
\end{corollary}
\begin{theorem}\label{thm33}
Let $\displaystyle G(\mathbb{Z}_{p^{n_{1}}_{1}p^{n_{2}}_{2}})$ be a unit graph, where $\displaystyle p_{1}$ and $p_{2}$ are odd distinct primes. Then
\begin{enumerate}
  \item $\displaystyle |V|=p^{n_{1}}_{1}p^{n_{2}}_{2}$ and $\displaystyle |E|=\frac{(p^{n_{1}}_{1}p^{n_{2}}_{2}-1)\phi(p^{n_{1}}_{1}p^{n_{2}}_{2})}{2}$.
  \item $\displaystyle \lambda(G(\mathbb{Z}_{p^{n_{1}}_{1}p^{n_{2}}_{2}}))=\delta(G(\mathbb{Z}_{p^{n_{1}}_{1}p^{n_{2}}_{2}}))=\phi(p^{n_{1}}_{1}p^{n_{2}}_{2})-1$.
\end{enumerate}
\begin{proof}
\begin{enumerate}
  \item Here $\displaystyle |V|=p^{n_{1}}_{1}p^{n_{2}}_{2}$. Since $\displaystyle 2\in U(\mathbb{Z}_{p^{n_{1}}_{1}p^{n_{2}}_{2}})$, by Theorem (\ref{thm1}), $\displaystyle \text{deg}(x)=\phi(p^{n_{1}}_{1}p^{n_{2}}_{2})-1$ for all $\displaystyle x\in U(\mathbb{Z}_{p^{n_{1}}_{1}p^{n_{2}}_{2}})$ and $\displaystyle \text{deg}(x)=\phi(p^{n_{1}}_{1}p^{n_{2}}_{2})$ for all $\displaystyle x\in N_{U}(\mathbb{Z}_{p^{n_{1}}_{1}p^{n_{2}}_{2}})$. Thus,
      \begin{align*}
        |E| &= \frac{\sum_{x\in V}\text{deg}(x)}{2}\\
            &= \frac{\sum_{x\in U(\mathbb{Z}_{p^{n_{1}}_{1}p^{n_{2}}_{2}})}\text{deg}(x)+\sum_{x\in N_{U}(\mathbb{Z}_{p^{n_{1}}_{1}p^{n_{2}}_{2}})}\text{deg}(x)}{2}\\
            &= \frac{\phi(p^{n_{1}}_{1}p^{n_{2}}_{2})(\phi(p^{n_{1}}_{1}p^{n_{2}}_{2})-1)+ (p^{n_{1}}_{1}p^{n_{2}}_{2}-\phi(p^{n_{1}}_{1}p^{n_{2}}_{2}))\phi(p^{n_{1}}_{1}p^{n_{2}}_{2})}{2}\\
        |E| &= \frac{(p^{n_{1}}_{1}p^{n_{2}}_{2}-1)\phi(p^{n_{1}}_{1}p^{n_{2}}_{2})}{2}.
      \end{align*}
  \item It follows from Theorem (\ref{thm32}) and Theorem (\ref{thm3}).
\end{enumerate}
\end{proof}
\end{theorem}

\begin{theorem}\label{thm34}
	Let $\displaystyle G(\mathbb{Z}_{2^{m}p^{n_{1}}_{1}})=(V,E)$ be a unit graph, where $p_{1}$ is an odd prime. Then $\displaystyle G(\mathbb{Z}_{2^{m}p^{n_{1}}_{1}})$ is a connected graph and $\displaystyle \text{diam}(G(\mathbb{Z}_{2^{m}p^{n_{1}}_{1}}))\leq 3$.
	\begin{proof}
		Let $\displaystyle \mathbb{Z}_{2^{m}p^{n_{1}}_{1}}=U(\mathbb{Z}_{2^{m}p^{n_{1}}_{1}})\cup N_{U}(\mathbb{Z}_{2^{m}p^{n_{1}}_{1}})$.Then 
		\begin{align*} N_{U}(\mathbb{Z}_{2^{m}p^{n_{1}}_{1}})&=N_{2}\cup N_{p_{1}}\cup N_{2p_{1}}\\
			\text{where  }  N_{2}&=\{x\in \mathbb{Z}_{2^{m}p^{n_{1}}_{1}} \ | \ x=2\alpha \ \& \ p_{1}\nmid \alpha \}\\
			 N_{p_{1}}&=\{x\in \mathbb{Z}_{2^{m}p^{n_{1}}_{1}} \ | \ x=p_{1}\alpha \ \& \ 2\nmid \alpha \}\\
			  N_{2p_{1}}&=\{x\in \mathbb{Z}_{2^{m}p^{n_{1}}_{1}} \ | \ x=2\alpha p_{1}  \}
			  \end{align*}
Now for any $x,y\in \mathbb{Z}_{2^{m}p^{n_{1}}_{1}}$, we have the following cases:\\
{\bf Case I:} If $\displaystyle x,y\in U(\mathbb{Z}_{2^{m}p^{n_{1}}_{1}})$, then $\displaystyle [x,0]$ and $\displaystyle [0,y]$ are edges in $\displaystyle G(\mathbb{Z}_{2^{m}p^{n_{1}}_{1}})$. Hence $\displaystyle d(x,y)=1$.\\
{\bf Case II:} If $\displaystyle x\in U(\mathbb{Z}_{2^{m}p^{n_{1}}_{1}})$ and $\displaystyle y\in N_{U}(\mathbb{Z}_{2^{m}p^{n_{1}}_{1}})$, then we have the following possibilities for $y$,\\
		(a) If $y\in N_{2}$, then $y=2\alpha$ for some integer $\alpha$. Suppose $x$ is not adjacent to $y$ which gives $\displaystyle p_{1}\mid x+2\alpha$. Since $p_{1}\in N_{P_{1}}$, which implies $N_{p_{1}}\neq \phi$. Let $\beta p_{1}\in N_{p_{1}}$. Then $\displaystyle z=2\alpha +\beta p_{1}\in U(\mathbb{Z}_{2^{m}p^{n_{1}}_{1}})$. Then $[z,y]$ is an edge. Since $\displaystyle x,z\in U(\mathbb{Z}_{2^{m}p^{n_{1}}_{1}})$, $x$ and $y$ are adjacent to $0$. Hence, we get $\displaystyle [x,0], \ [0,z]$ and $\displaystyle [z,y]$ are edges in $\displaystyle G(\mathbb{Z}_{2^{m}p^{n_{1}}_{1}})$, which gives $d(x,y)\leq 3$.\\
		(b) If $\displaystyle y\in N_{p_{1}}$, then $y=\alpha p_{1}$. Note that $x$ is not adjacent to $y$. Since $\displaystyle 2\in N_{2}$ which implies $N_{2}\neq \phi$. Now for any, $z=2\beta\in N_{2}$,  $y$ is adjacent to $z$. Now it is enough to show that $x$ is adjacent to $\displaystyle 2\beta$ or $\displaystyle 2\beta$, suppose $x$ is not adjacent to $\displaystyle 2\beta \in N_{2}$ and $\displaystyle  -2\beta \in N_{2}$, which implies $\displaystyle p_{1}\mid x+2\beta$ and $\displaystyle p_{1}\mid x-2\beta$ which gives $p_{1}\mid 4\beta$,  hence $p_{1}\mid \beta$, which is contradiction to choice of $\displaystyle \beta$. Hence,  $\displaystyle [x,2\beta]$ or $\displaystyle [x,-2\beta]$ is edges in $\displaystyle G(\mathbb{Z}_{2^{m}p^{n_{1}}_{1}})$, which implies $d(x,y)=2$.\\
		(c) If $\displaystyle y\in N_{2p_{1}}$, then $y=2\alpha p_{1}$. Since $\displaystyle x\in U(\mathbb{Z}_{2^{m}p^{n_{1}}_{1}})$, $x$ is adjacent to $y$, if this not hold then either $\displaystyle p_{1}\mid x+y$ or $\displaystyle 2\mid x+y$. If $\displaystyle p_{1}\mid x+2\alpha p_{1}$, then $p_{1}\mid x$ which is a contradiction. In similar way, If $\displaystyle 2\mid x+2\alpha p_{1}$, then $\displaystyle 2\mid x$ which is also contradiction to $\displaystyle x\in U(\mathbb{Z}_{2^{m}p^{n_{1}}_{1}})$. Hence $\displaystyle [x,y]$ is an edge in $\displaystyle G(\mathbb{Z}_{2^{m}p^{n_{1}}_{1}})$ which gives, $\displaystyle d(x,y)=1\leq 3$.\\
		{\bf Case III:} If $\displaystyle x,y\in N_{U}(\mathbb{Z}_{2^{m}p^{n_{1}}_{1}})$, then we have the following possibilities for $x$ and $y$:\\
		(a) If $\displaystyle x,y\in N_{2}$, then $\displaystyle x=2\alpha$ and $\displaystyle y=2\alpha'$. For any $\displaystyle z=\beta p_{1}\in N_{p_{1}}$, we have $[x,z]$ and $[z,y]$ are edges in $\displaystyle G(\mathbb{Z}_{2^{m}p^{n_{1}}_{1}})$. Hence, we have $ \displaystyle d(x,y)=2$.\\
		(b) If $\displaystyle x,y\in N_{p_{1}}$, then from above procedure for any $\displaystyle z\in N_{2}$, we have $[x,z]$ and $[z,y]$ are edges in $\displaystyle G(\mathbb{Z}_{2^{m}p^{n_{1}}_{1}})$. Hence we get $\displaystyle d(x,y)=2\leq 3$.\\
		(c) If $\displaystyle x,y\in N_{2p_{1}}$, then for any $\displaystyle z\in U(\mathbb{Z}_{2^{m}p^{n_{1}}_{1}})$, we have $\displaystyle [x,z]$ and $\displaystyle [z,y]$ are edges in $\displaystyle G(\mathbb{Z}_{2^{m}p^{n_{1}}_{1}})$, which gives $\displaystyle d(x,y)=2$.\\
		(d) If $\displaystyle x\in N_{2}$ and $\displaystyle y\in N_{p_{1}}$ or vice versa, then $x$ is adjacent to $y$, which implies $d(x,y)=1$.\\
		(e)If $\displaystyle x\in N_{2}$ and $y\in N_{2p_{1}}$, then $x=2\alpha$ and $\displaystyle y=2\alpha' p_{1} $. Now for any $z=\beta p_{1}\in N_{p_{1}}$, we have $\displaystyle x+z\in N_{U}(\mathbb{Z}_{2^{m}p^{n_{1}}_{1}})$ and $x$ is adjacent to $x+z$. Since $\displaystyle x+z\in U(\mathbb{Z}_{2^{m}p^{n_{1}}_{1}})$, we have $x+z$ is adjacent to $y$. Hence, we get $[x,x+z]$ and $[x+z,y]$ are edges in $\displaystyle G(\mathbb{Z}_{2^{m}p^{n_{1}}_{1}})$, which gives $d(x,y)=2\leq 3$.\\
		(f) If $\displaystyle x\in N_{p_{1}}$ and $y\in N_{2p_{1}}$, then $x=\alpha p_{1}$ and $\displaystyle y=2\alpha' p_{1} $. Now for any $z=2\beta \in N_{2}$, we have $x$ is adjacent to $z$ and $z$ is adjacent to $x+z$ both follows from above part. Note that $\displaystyle x+z\in U(\mathbb{Z}_{2^{m}p^{n_{1}}_{1}})$. Hence, we have $x+z$ is adjacent to $y$, from this, we get $\displaystyle [x,z], \ [z,x+z]$ and $\displaystyle [x+z,y]$ are edges in $\displaystyle G(\mathbb{Z}_{2^{m}p^{n_{1}}_{1}})$, which give $d(x,y)\leq 3$.\\
		Hence, $\displaystyle G(\mathbb{Z}_{2^{m}p^{n_{1}}_{1}})$ is connected graph and $\displaystyle d(x,y)\leq 3$, for all $\displaystyle x,y\in \mathbb{Z}_{2^{m}p^{n_{1}}_{1}}$. Therefore, we have  $\displaystyle \text{diam}(G(\mathbb{Z}_{2^{m}p^{n_{1}}_{1}}))\leq 3$.
	\end{proof}
\end{theorem}

\begin{corollary}\label{thm37}
Let $\displaystyle G(\mathbb{Z}_{2^{m}p^{n_{1}}_{1}})$ be a unit graph, where $p_{1}$ is an odd prime.
\begin{enumerate}
  \item If $\displaystyle 2^{m}p^{n_{1}}_{1}\neq 6$ then $\displaystyle g_{r}(G(\mathbb{Z}_{2^{m}p^{n_{1}}_{1}}))=4$
  \item If $\displaystyle 2^{m}p^{n_{1}}_{1}= 6$ then $\displaystyle g_{r}(G(\mathbb{Z}_{2^{m}p^{n_{1}}_{1}}))=6$
\end{enumerate}
\begin{proof}
Proof follows from Theorems (\ref{thm7}) and (\ref{thm31}).
\end{proof}
\end{corollary}

\begin{theorem}\label{thm35}
	Let $\displaystyle G(\mathbb{Z}_{2^{m}p^{n_{1}}_{1}})$ be a unit graph, where $\displaystyle p_{1}$ is an odd prime. Then
	\begin{enumerate}
		\item $\displaystyle |V|=2^{m}p^{n_{1}}_{1}$ and $\displaystyle |E|=2^{m-1}p^{n_{1}}_{1}\phi(2^{m}p^{n_{1}}_{1})$.
		\item  $\displaystyle \lambda(G(\mathbb{Z}_{2^{m}p^{n_{1}}_{1}}))=\delta(G(\mathbb{Z}_{2^{m}p^{n_{1}}_{1}}))=\phi(2^{m}p^{n_{1}}_{1})$.
	\end{enumerate}
	\begin{proof}
		\begin{enumerate}
			\item From the definition of unit graph, we have $\displaystyle |V|=2^{m}p^{n_{1}}_{1}$. Since $\displaystyle 2\in N_{U}(\mathbb{Z}_{2^{m}p^{n_{1}}_{1}})$ from Theorem (\ref{thm1}), we get $\displaystyle G(\mathbb{Z}_{2^{m}p^{n_{1}}_{1}})$ is $\displaystyle \phi(2^{m}p^{n_{1}}_{1})$-regular graph and hence, we have $\displaystyle |E|=\frac{2^{m}p^{n_{1}}_{1}\phi(2^{m}p^{n_{1}}_{1})}{2}=2^{m-1}p^{n_{1}}_{1}\phi(2^{m}p^{n_{1}}_{1})$.
			\item From Lemma (\ref{thm31}) and Theorems (\ref{thm33}) \& (\ref{thm4}), we get $\displaystyle \lambda(G(\mathbb{Z}_{2^{m}p^{n_{1}}_{1}}))=\delta(G(\mathbb{Z}_{2^{m}p^{n_{1}}_{1}}))=\phi(2^{m}p^{n_{1}}_{1})$.
		\end{enumerate}
	\end{proof}
\end{theorem}

\begin{theorem} Let $\displaystyle G(\mathbb{Z}_{n})$ be a unit graph, where $n=p^{n_{1}}_{1}p^{n_{2}}_{2}$ and $p_{1} \ \& \ p_{2}$ are distinct primes. Let $H$ be a $|V|\times |E|$ incidence matrix of $G(\mathbb{Z}_{n})$.
    \begin{enumerate}
      \item If $\displaystyle 2\in U(\mathbb{Z}_{n})$, then binary code generated by $H$ is a  $\displaystyle C_{2}(H)=\left[\frac{(n-1)\phi(n)}{2},n-1,\phi(n)-1\right]_{2}$ code over finite field $\mathbb{F}_{2}$.
      \item If $\displaystyle 2\in N_{U}(\mathbb{Z}_{n})$, then for any odd prime $q$, the $q$-ary code generated by $H$ is a $\displaystyle C_{q}(H)=\left[\frac{n\phi(n)}{2},n-1,\phi(n)\right]_{q}$ code over finite field $\mathbb{F}_{q}$.
    \end{enumerate}
\begin{proof}
\begin{enumerate}
  \item If $\displaystyle 2\in U(\mathbb{Z}_{p^{n_{1}}_{1}p^{n_{2}}_{2}})$, then $p_{1}$ and $p_{2}$ are odd primes. By Theorem (\ref{thm32}), $\displaystyle G(\mathbb{Z}_{p^{n_{1}}_{1}p^{n_{2}}_{2}})=(V,E)$ is a connected graph and hence by Theorem (\ref{thm2}), binary code generated by $H$ is
      $\displaystyle C_{2}(H)=[|E|,|V|-1,\lambda(G(\mathbb{Z}_{p^{n_{1}}_{1}p^{n_{2}}_{2}}))]_{2}$. Now from Theorem (\ref{thm33}), we get $\displaystyle |E|= \frac{(n-1)\phi(n)}{2},\\  |V|-1=n-1$ and $\displaystyle \lambda(G(\mathbb{Z}_{p^{n_{1}}_{1}p^{n_{2}}_{2}}))=\phi(n)-1$.
  \item If $\displaystyle 2\in N_{U}(\mathbb{Z}_{p^{n_{1}}_{1}p^{n_{2}}_{2}})$, then either $p_{1}=2$ or $p_{2}=2$. By Theorem (\ref{thm34}) and Lemma (\ref{thm31}), $\displaystyle G(\mathbb{Z}_{p^{n_{1}}_{1}p^{n_{2}}_{2}})=(V,E)$ is a connected bipartite graph and hence by Theorems (\ref{thm2}) and (\ref{thm5}), $q$-ary code generated by $H$ is $\displaystyle C_{q}(H)=[|E|,|V|-1,\lambda(G(\mathbb{Z}_{2^{m}p^{n_{1}}_{1}}))]_{q}$. Now using Theorem (\ref{thm35}), we conclude the result.
\end{enumerate}
\end{proof}
\end{theorem}

\begin{corollary}
Let $C_{q}(H)$ and $C_{2}(H)$ denote the linear codes generated from incidence matrices of $\displaystyle G(\mathbb{Z}_{2^{m}p^{n_{1}}_{1}})$ and $\displaystyle G(\mathbb{Z}_{p^{n_{1}}_{1}p^{n_{2}}_{2}})$. Then
\begin{enumerate}
  \item Dual of code $C_{2}$ is $\displaystyle C^{\perp}_{2}=\left[\frac{(n-1)\phi(n)}{2},\frac{(n-1)[\phi(n)-2]}{2},3\right]_{2}$, where $\displaystyle n=p^{n_{1}}_{1}p^{n_{2}}_{2}$.
  \item Dual of code $C_{q}$ is $\displaystyle C^{\perp}_{q}=\left[\frac{n\phi(n)}{2},\frac{n(\phi(n)-2)+2}{2},4\right]_{q}$, where $\displaystyle n=2^{m}p^{n_{1}}_{1}\neq 6$.
\end{enumerate}
\begin{proof}
\begin{enumerate}
\item From Theorem (\ref{thm8}), $\text{diam}(\displaystyle C^{\perp}_{2})=\frac{(n-1)[\phi(n)-2]}{2}$.\\
 By Theorem (\ref{thm6}) and Corollary  (\ref{thm36}), $\displaystyle d(C^{\perp}_{2})=3$.
  \item Proof follows from Theorems (\ref{thm8}) \& (\ref{thm6}) and Corollary (\ref{thm37}).
\end{enumerate}
\end{proof}
\end{corollary}

\section{Construction of Linear Codes from the Incidence Matrix of a Unit Graph $G(\mathbb{Z}_{p_{1}^{n_{1}}p_{2}^{n_{2}}p_{3}^{n_{3}}})$}
In this section, we extend the result for three distinct primes $p_{1}$, $p_{2}$ and $p_{3}$.
\begin{theorem}\label{thm41}
	Let $\displaystyle G(\mathbb{Z}_{p^{n_{1}}_{1}p^{n_{2}}_{2}p^{n_{3}}_{3}})$ be a unit graph, where $\displaystyle p_{1}, \ p_{2}$ and  $\displaystyle p_{3}$ are distinct odd primes. Then $\displaystyle G(\mathbb{Z}_{p^{n_{1}}_{1}p^{n_{2}}_{2}p^{n_{3}}_{3}})$ is connected graph and $\displaystyle \text{diam}(G(\mathbb{Z}_{p^{n_{1}}_{1}p^{n_{2}}_{2}p^{n_{3}}_{3}}))\leq 2$.
	\begin{proof}
		Let $\mathbb{Z}_{p^{n_{1}}_{1}p^{n_{2}}_{2}p^{n_{3}}_{3}}=U(\mathbb{Z}_{p^{n_{1}}_{1}p^{n_{2}}_{2}p^{n_{3}}_{3}})\cup N_{U}(\mathbb{Z}_{p^{n_{1}}_{1}p^{n_{2}}_{2}p^{n_{3}}_{3}})$. Also we can express $N_{U}(\mathbb{Z}_{p^{n_{1}}_{1}p^{n_{2}}_{2}p^{n_{3}}_{3}})=\displaystyle \bigcup_{i=1}^{3}N_{p_{i}}\bigcup_{i\neq j}N_{p_{i}p_{j}}\bigcup N_{p_{1}p_{2}p_{3}}$,\\ where
		$ \displaystyle N_{p_{i}}=\left\{x\in \mathbb{Z}_{p^{n_{1}}_{1}p^{n_{2}}_{2}p^{n_{3}}_{3}} \ | \ x=\alpha p_{i} \ \& \ p_{j}\nmid \alpha, \ \forall \ i\neq j \right\}$,\\ $\displaystyle N_{p_{i}p_{j}}=\left\{x\in \mathbb{Z}_{p^{n_{1}}_{1}p^{n_{2}}_{2}p^{n_{3}}_{3}} \ | \ x=\alpha p_{i}p_{j} \ \& \ p_{k}\nmid \alpha, \ \ k\neq i,j\right\}$\\ and  $\displaystyle N_{p_{1}p_{2}p_{3}}=\left\{x\in \mathbb{Z}_{p^{n_{1}}_{1}p^{n_{2}}_{2}p^{n_{3}}_{3}} \ | \ x=\alpha p_{1}p_{2}p_{3}\right\}$. Now for any $\displaystyle x,y\in \mathbb{Z}_{p^{n_{1}}_{1}p^{n_{2}}_{2}p^{n_{3}}_{3}}$ We have the following cases:\\
		{\bf Case I:} If $\displaystyle x,y\in U(\mathbb{Z}_{p^{n_{1}}_{1}p^{n_{2}}_{2}p^{n_{3}}_{3}})$, then $\displaystyle [x,0]$ and $\displaystyle [0,y]$ are edges. Hence, $d(x,y)\leq 2$.\\
		{\bf Case II:} If $\displaystyle x\in U(\mathbb{Z}_{p^{n_{1}}_{1}p^{n_{2}}_{2}p^{n_{3}}_{3}})$ and $\displaystyle y\in N_{U}(\mathbb{Z}_{p^{n_{1}}_{1}p^{n_{2}}_{2}p^{n_{3}}_{3}})$, then we have the following possibilities for $y$:\\
		(a) If $\displaystyle y\in N_{p_{i}}$, then $y=\alpha p_{i}$. Suppose $x$ and $y$ are not adjacent then $\displaystyle p_{j}\mid x+\alpha p_{i}$ for $i\neq j$ which implies $ x+\alpha p_{i}=\beta p_{j}$. Note that $\displaystyle p_{i}\nmid \beta$, hence $\beta p_{j}p_{k}\in N_{p_{j}p_{k}}$ and $-\beta p_{j}p_{k}\in N_{p_{j}p_{k}}$  for $i\neq j,k$. Now our claim is either $\displaystyle [x,\beta p_{j}p_{k}]$ or $[x,-\beta p_{j}p_{k}]$ is edge in $G(\mathbb{Z}_{p^{n_{1}}_{1}p^{n_{2}}_{2}p^{n_{3}}_{3}})$. Suppose both not hold which implies $\displaystyle p_{i}\mid x+\beta p_{j}p_{k}$ and $\displaystyle p_{i}\mid x-\beta p_{j}p_{k}$, from this we get  $\displaystyle p_{i}\mid 2x$, this implies $p_{i}\mid x$ which is a contradiction.  Now we prove $y$ is adjacent to $\beta p_{j}p_{k}$, suppose this is not true, which implies $\displaystyle p_{i}\mid \alpha p_{i}+\beta p_{j}p_{k}$,  from this we have, $p_{i}\mid \beta$ which is a contradiction. Similar contradiction occurs, when we assuming $\displaystyle p_{j}\mid \alpha p_{i}+\beta p_{j}p_{k}$ and $\displaystyle p_{k}\mid \alpha p_{i}+\beta p_{j}p_{k}$. Hence $\displaystyle [\alpha p_{i},\beta p_{j}p_{k}]$ is an edge. In similar way, we can show that $\displaystyle [\alpha p_{i},-\beta p_{j}p_{k}]$ is an edge in $\displaystyle G(\mathbb{Z}_{p^{n_{1}}_{1}p^{n_{2}}_{2}p^{n_{3}}_{3}})$. Hence, there exist $\displaystyle  z\in N_{p_{j}p_{k}}$ such that $\displaystyle [x,z]$ and $\displaystyle [z,y]$ are edges in $G(\mathbb{Z}_{p^{n_{1}}_{1}p^{n_{2}}_{2}p^{n_{3}}_{3}})$, from this we concluded that $\displaystyle d(x,y)\leq 2$.\\
		(b) If $\displaystyle y\in N_{p_{i}p_{j}}$, then $y=\alpha p_{i}p_{j}$. Suppose $x$ and $y$ are not adjacent then $\displaystyle p_{k}\mid x+\alpha p_{i}p_{j}$ for $k\neq i,j$, which implies $ x+\alpha p_{i}p_{j}=\beta p_{k}$. Note that $\displaystyle p_{i}\nmid \beta$ and $\displaystyle p_{j}\nmid \beta$ which gives $\beta p_{k}\in N_{p_{k}}$.  Now we prove that $\displaystyle [x,\beta p_{k}]$ is an edge, suppose this is not an edge, which implies $\displaystyle p_{i}\mid x+\beta p_{k}$ for some $i\neq k$. Hence $\displaystyle p_{i}\mid x+x+\alpha p_{i}p_{j}$ this implies $p_{i}\mid 2x$, which is a contradiction. $\displaystyle [x,\beta p_{k}]$ is an edge and using procedure in $(a)$, we can show that $\displaystyle [\beta p_{k}, \alpha p_{i}p_{j}]=[\beta p_{k},y]$. Hence $d(x,y)\leq 2$.\\
		(c) If $\displaystyle y\in N_{p_{1}p_{2}p_{3}}$, then $\displaystyle [x,\alpha p_{1}p_{2}p_{3}]$ is an edge in $G(\mathbb{Z}_{p^{n_{1}}_{1}p^{n_{2}}_{2}p^{n_{3}}_{3}})$. Hence $d(x,y)=1$.\\
		{\bf Case III:} If $\displaystyle x,y\in N_{U}(\mathbb{Z}_{p^{n_{1}}_{1}p^{n_{2}}_{2}p^{n_{3}}_{3}})$, then we have the following possibilities:\\
		(a) If $\displaystyle x,y\in N_{p_{i}}$, then $x=\alpha p_{i}$ and $\displaystyle y=\alpha' p_{i}$. Since $N_{p_{j}p_{k}}\neq \phi$. For any $\displaystyle \beta p_{j}p_{k}\in N_{p_{j}p_{k}}$ for $i\neq j,k$, we have $\displaystyle [x,\beta p_{j} p_{k}]$ and $\displaystyle [\beta p_{j} p_{k},y]$ are edges in $G(\mathbb{Z}_{p^{n_{1}}_{1}p^{n_{2}}_{2}p^{n_{3}}_{3}})$. Hence, $d(x,y)\leq 2$.\\
		(b) If $\displaystyle x,y\in N_{p_{i}p_{j}}$, then $x=\alpha p_{i}p_{j}$ and $\displaystyle y=\alpha' p_{i}p_{j}$. Let $\displaystyle \beta p_{k}\in N_{p_{k}}$ for $k\neq i,j$. Then $\displaystyle [x,\beta p_{k}]$ and $\displaystyle [\beta p_{k},y]$ are edges in $G(\mathbb{Z}_{p^{n_{1}}_{1}p^{n_{2}}_{2}p^{n_{3}}_{3}})$. Hence, $d(x,y)\leq 2$.\\
		(c) $\displaystyle x,y\in N_{p_{1}p_{2}p_{3}}$, then for any $\displaystyle z\in U(\mathbb{Z}_{p^{n_{1}}_{1}p^{n_{2}}_{2}p^{n_{3}}_{3}})$, $\displaystyle [x,z]$ and $\displaystyle [z,x]$ are edges in $G(\mathbb{Z}_{p^{n_{1}}_{1}p^{n_{2}}_{2}p^{n_{3}}_{3}})$. Hence, $d(x,y)\leq 2$.\\
		(d) If $\displaystyle x\in N_{p_{i}}$ and $\displaystyle y\in N_{p_{j}}$ for $i\neq j$, then $\displaystyle x=\alpha p_{i}$ and $y=\alpha' p_{j}$. Assume that $\displaystyle [x,y]$ is not an edge, which gives $\displaystyle p_{k}\mid x+y$, for $k\neq i,j$. Now our claim is that $\displaystyle [x+y,x]$ and $\displaystyle [y,x+y]$ are edges. Suppose $x+y$ is not adajcent to $x$, hence $ p_{m}\mid 2x+y$ for some $\displaystyle m=1,2,3$.  If $ p_{i}\mid 2x+y$ i.e. $p_{i}\mid 2\alpha p_{i}+\alpha' p_{j}$,  then $p_{i}\mid \alpha'$, which is a contradiction. Similarly, if $ p_{j}\mid 2x+y$,  then $p_{j}\mid \alpha$ and if $ p_{k}\mid 2x+y$ i.e. $p_{k}\mid x+y+x$, then $p_{k}\mid x$, since $\displaystyle p_{k}\mid x+y$. Hence, $\displaystyle [x+y,x]$ is an edge in $G(\mathbb{Z}_{p^{n_{1}}_{1}p^{n_{2}}_{2}p^{n_{3}}_{3}})$. In similar way, we can show that $x+y$ is adjacent to $y$, which gives $\displaystyle d(x,y)\leq 2$.\\
		(e) If $\displaystyle x\in N_{p_{i}p_{j}}$ and $\displaystyle y\in N_{p_{j}}$, then $\displaystyle x=\alpha p_{i}p_{j}$ and $\displaystyle y=\alpha' p_{j}$. Note that for any $\displaystyle z=\beta p_{k}$, for $k\neq i,j$, $x$ is adjacent to $z$. Now assuming that $y=\alpha p_{j}\in N_{p_{j}}$ and $\displaystyle \beta p_{k}\in N_{p_{k}}$ are not adjacent, then $\displaystyle p_{i}\mid \alpha p_{j}+\beta p_{k}$. Since $\displaystyle \beta p_{k}\in N_{p_{k}}$, which implies $-\beta p_{k}\in N_{p_{k}}$ this gives $\displaystyle [y,-\beta p_{k}]$ is an edge. Hence, $y$ is adjacent to either $\displaystyle \beta p_{k}$ or $\displaystyle -\beta p_{k}$, which gives $d(x,y)\leq 2$.\\
		(f) If $\displaystyle x\in N_{p_{i}p_{j}}$ and $\displaystyle y\in N_{p_{j}p_{k}}$ for $\displaystyle i\neq k$, then $\displaystyle  x=\alpha p_{i}p_{j}$ and $\displaystyle  y=\alpha'p_{j}p_{k}$. Now for any $\displaystyle \beta p_{i}\in N_{p_{i}}$, we have $\displaystyle z=(\beta p_{i}+\alpha' p_{j})p_{k}+\alpha p_{i}p_{j}\in U(\mathbb{Z}_{p^{n_{1}}_{1}p^{n_{2}}_{2}p^{n_{3}}_{3}})$ such that $\displaystyle [x,z]$ and $\displaystyle [z,y]$ are edges in $G(\mathbb{Z}_{p^{n_{1}}_{1}p^{n_{2}}_{2}p^{n_{3}}_{3}})$. To see this, suppose $x$ is not adjacent to $z$ which implies $\displaystyle p_{i}\mid \alpha p_{i}p_{j}+(\beta p_{i}+\alpha' p_{j})p_{k}+\alpha p_{i}p_{j}$, we get $p_{i}\mid \alpha'$ which is not possible. Similarly , if $\displaystyle p_{j}\mid \alpha p_{i}p_{j}+(\beta p_{i}+\alpha' p_{j})p_{k}+\alpha p_{i}p_{j}$, then $p_{j}\mid \beta$ and if $\displaystyle p_{k}\mid \alpha p_{i}p_{j}+(\beta p_{i}+\alpha' p_{j})p_{k}+\alpha p_{i}p_{j}$, then $p_{k}\mid 2\alpha$ both are not possible, which gives $x$ is adjacent to $z$.  Hence, we get $d(x,y)\leq 2$.\\
		(g) If $\displaystyle x\in N_{p_{i}p_{j}}$ and $\displaystyle y\in N_{p_{k}}$ for $k\neq i,j$, then $\displaystyle [x,y]$ is an edge in $G(\mathbb{Z}_{p^{n_{1}}_{1}p^{n_{2}}_{2}p^{n_{3}}_{3}})$ gives $\displaystyle d(x,y)=1$.\\
Hence, $\displaystyle G(\mathbb{Z}_{p^{n_{1}}_{1}p^{n_{2}}_{2}p^{n_{3}}_{3}})$ is connected graph and $d(x,y)\leq 2$, for all $\displaystyle x,y\in \mathbb{Z}_{p^{n_{1}}_{1}p^{n_{2}}_{2}p^{n_{3}}_{3}}$, which gives, $\displaystyle \text{diam}(G(\mathbb{Z}_{p^{n_{1}}_{1}p^{n_{2}}_{2}p^{n_{3}}_{3}}))\leq 2$.
	\end{proof}
\end{theorem}
\begin{corollary}\label{thm45}
	Let $\displaystyle G(\mathbb{Z}_{p^{n_{1}}_{1}p^{n_{2}}_{2}p^{n_{3}}_{3}})$ be a unit graph, where $\displaystyle p_{1},p_{2}$ and  $\displaystyle p_{3}$ are odd primes. Then $\displaystyle g_{r}(G(\mathbb{Z}_{p^{n_{1}}_{1}p^{n_{2}}_{2}p^{n_{3}}_{3}}))=3$.
	\begin{proof}
The result follows from $p_{1}$ is adjacent to $p_{2}p_{3}$ and $p_{1}+p_{2}p_{3}$ is adjacent to  both $p_{1} \ \& \ p_{2}p_{3}$.
	\end{proof}
\end{corollary}
\begin{theorem}\label{thm42}
	Let $\displaystyle G(\mathbb{Z}_{p^{n_{1}}_{1}p^{n_{2}}_{2}p^{n_{3}}_{3}})$ be a unit graph, where $\displaystyle p_{1},p_{2}$ and $p_{3}$ are odd primes. Then
	\begin{enumerate}
		\item $\displaystyle |V|=p^{n_{1}}_{1}p^{n_{2}}_{2}p^{n_{3}}_{3}$ and $\displaystyle |E|=\frac{(p^{n_{1}}_{1}p^{n_{2}}_{2}p^{n_{3}}_{3}-1)\phi(p^{n_{1}}_{1}p^{n_{2}}_{2}p^{n_{3}}_{3})}{2}$.
		\item $\displaystyle \lambda(G(\mathbb{Z}_{p^{n_{1}}_{1}p^{n_{2}}_{2}p^{n_{3}}_{3}}))=\delta(G(\mathbb{Z}_{p^{n_{1}}_{1}p^{n_{2}}_{2}p^{n_{3}}_{3}}))=\phi(p^{n_{1}}_{1}p^{n_{2}}_{2}p^{n_{3}}_{3})-1$.
	\end{enumerate}
	\begin{proof}
		\begin{enumerate}
			\item Clearly, $\displaystyle |V|=p^{n_{1}}_{1}p^{n_{2}}_{2}p^{n_{3}}_{3}$. Since $\displaystyle 2\in U(\mathbb{Z}_{p^{n_{1}}_{1}p^{n_{2}}_{2}p^{n_{3}}_{3}})$ by Theorem (\ref{thm1}) $\displaystyle \text{deg}(x)=\phi(p^{n_{1}}_{1}p^{n_{2}}_{2}p^{n_{3}}_{3})-1$ for all $\displaystyle x\in U(\mathbb{Z}_{p^{n_{1}}_{1}p^{n_{2}}_{2}p^{n_{3}}_{3}})$ and $\displaystyle \text{deg}(x)=\phi(p^{n_{1}}_{1}p^{n_{2}}_{2}p^{n_{3}}_{3})$ for all $\displaystyle x\in N_{U}(\mathbb{Z}_{p^{n_{1}}_{1}p^{n_{2}}_{2}p^{n_{3}}_{3}})$. Thus,
			\begin{align*}
				|E| &= \frac{\sum_{x\in V}\text{deg}(x)}{2}\\
				&= \frac{\sum_{x\in U(\mathbb{Z}_{p^{n_{1}}_{1}p^{n_{2}}_{2}p^{n_{3}}_{3}})}\text{deg}(x)+\sum_{x\in N_{U}(\mathbb{Z}_{p^{n_{1}}_{1}p^{n_{2}}_{2}p^{n_{3}}_{3}})}\text{deg}(x)}{2}\\
				&= \frac{\phi(p^{n_{1}}_{1}p^{n_{2}}_{2}p^{n_{3}}_{3})(\phi(p^{n_{1}}_{1}p^{n_{2}}_{2}p^{n_{3}}_{3})-1)+ (p^{n_{1}}_{1}p^{n_{2}}_{2}p^{n_{3}}_{3}-\phi(p^{n_{1}}_{1}p^{n_{2}}_{2}p^{n_{3}}_{3}))\phi(p^{n_{1}}_{1}p^{n_{2}}_{2}p^{n_{3}}_{3})}{2}\\
				|E| &= \frac{(p^{n_{1}}_{1}p^{n_{2}}_{2}p^{n_{3}}_{3}-1)\phi(p^{n_{1}}_{1}p^{n_{2}}_{2}p^{n_{3}}_{3})}{2}.
			\end{align*}
			\item It follows from Theorem (\ref{thm41}) and Theorem (\ref{thm3}).
		\end{enumerate}
	\end{proof}
\end{theorem}
\begin{theorem}\label{thm43}
Let $\displaystyle G(\mathbb{Z}_{2^{m}p^{n_{1}}_{1}p^{n_{2}}_{2}})$ be a unit graph, where $p_{1}$ and $p_{2}$ are odd primes. Then $\displaystyle G(\mathbb{Z}_{2^{m}p^{n_{1}}_{1}p^{n_{2}}_{2}})$ is connected graph and $\displaystyle \text{diam}( G(\mathbb{Z}_{2^{m}p^{n_{1}}_{1}p^{n_{2}}_{2}}))\leq 3$.
\begin{proof}
Let $\displaystyle \mathbb{Z}_{2^{m}p^{n_{1}}_{1}p^{n_{2}}_{2}}=U(\mathbb{Z}_{2^{m}p^{n_{1}}_{1}p^{n_{2}}_{2}})\cup N_{U}(\mathbb{Z}_{2^{m}p^{n_{1}}_{1}p^{n_{2}}_{2}})$ and\\ $\displaystyle N_{U}(\mathbb{Z}_{2^{m}p^{n_{1}}_{1}p^{n_{2}}_{2}})=N_{2}\cup N_{p_{1}}\cup N_{p_{2}}\cup N_{2p_{1}}\cup N_{2p_{2}}\cup N_{p_{1}p_{2}}\cup N_{2p_{1}p_{2}}$\\
 where $\displaystyle N_{2}=\{x\in \mathbb{Z}_{2^{m}p^{n_{1}}_{1}p^{n_{2}}_{2}} \ | \ x=2\alpha \ \& \ p_{1}\nmid \alpha \ \& \ p_{2}\nmid \alpha \}$\\
  $\displaystyle N_{p_{i}}=\{x\in \mathbb{Z}_{2^{m}p^{n_{1}}_{1}p^{n_{2}}_{2}} \ | \ x=\alpha p_{i} \ \& \ 2\nmid \alpha \ \& \ p_{j}\nmid \alpha, \ i\neq j  \}$\\ $\displaystyle N_{2p_{i}}=\{x\in \mathbb{Z}_{2^{m}p^{n_{1}}_{1}p^{n_{2}}_{2}} \ | \ x=2\alpha p_{i} \ \& \ p_{j}\nmid \alpha, \ i\neq j \}\\
   N_{p_{1}p_{2}}=\{x\in \mathbb{Z}_{2^{m}p^{n_{1}}_{1}p^{n_{2}}_{2}} \ |  \ x=\alpha p_{1}p_{2} \ \& \ 2\nmid \alpha\}$\\
    and $N_{2p_{1}p_{2}}=\{x\in \mathbb{Z}_{2^{m}p^{n_{1}}_{1}p^{n_{2}}_{2}} \ |  \ x=\alpha 2p_{1}p_{2}\}$.\\
  Now for any $x,y\in \mathbb{Z}_{2^{m}p^{n_{1}}_{1}p^{n_{2}}_{2}}$, we have the following cases:\\
{\bf Case I:} If $\displaystyle x,y\in U(\mathbb{Z}_{2^{m}p^{n_{1}}_{1}p^{n_{2}}_{2}})$, then $x$ and $y$ adjacent to $0$, which gives $\displaystyle d(x,y)=1\leq 3$.\\
{\bf Case II:} If  $\displaystyle x\in U(\mathbb{Z}_{2^{m}p^{n_{1}}_{1}p^{n_{2}}_{2}})$ and $\displaystyle y\in N_{U}(\mathbb{Z}_{2^{m}p^{n_{1}}_{1}p^{n_{2}}_{2}})$, then we have the following possibilities for $y$:\\
(a) If $\displaystyle y\in N_{2}$, then $\displaystyle y=2\alpha$. Suppose $x$ is not adjacent to $y$. For any $\displaystyle z=\beta p_{1}p_{2}$, we have $\displaystyle y+z \in U(\mathbb{Z}_{2^{m}p^{n_{1}}_{1}p^{n_{2}}_{2}})$. Now our claim is that, $\displaystyle [y,y+z]$ is an edge in $\displaystyle G(\mathbb{Z}_{2^{m}p^{n_{1}}_{1}p^{n_{2}}_{2}})$. If this is not true, then $\displaystyle 2\mid 2y+z$, which gives $2\mid \beta$. Similarly, if $\displaystyle p_{i}\mid 4\alpha+\beta p_{1}p_{2}$, then $p_{i}\mid \alpha$ for $i=1,2$ in both cases we get contradiction. Since,  $\displaystyle y+z,x \in U(\mathbb{Z}_{2^{m}p^{n_{1}}_{1}p^{n_{2}}_{2}})$, which gives $[y+z,0]$ and $[0,x]$ are edges in $\displaystyle G(\mathbb{Z}_{2^{m}p^{n_{1}}_{1}p^{n_{2}}_{2}})$. Hence, we get $d(x,y)\leq 3$.\\
(b) If $y\in N_{p_{i}}$, for $i=1,2$, then $x$ is not adjacent to $y$. First we show that $x$ is adjacent to some $z =2\beta p_{j}\in N_{2p_{j}}$ for $i\neq j$. Suppose $x$ is not adjacent to $z=2\beta p_{j}$, which implies $\displaystyle p_{i}\mid x+z$. Note that
      $\displaystyle p_{j}\nmid x+z$ and $\displaystyle 2\nmid x+z$. Since $\displaystyle z=2\beta p_{j}\in N_{2p_{j}}$, which gives $ -2\beta p_{j}\in N_{2p_{j}}$. If $\displaystyle p_{i}\mid x-2\beta p_{j}$, then $\displaystyle p_{i}\mid \beta$, which is not possible. Hence, there exits $\displaystyle z\in N_{2p_{j}}$ such that $x$ is adjacent to $z$. Clearly, $y$ is adjacent to $z$, for all $\displaystyle z\in N_{2p_{j}}$. Hence, we get $\displaystyle d(x,y)=2$.\\
(c) If $\displaystyle y\in N_{2p_{i}}$, then $ y=2\alpha p_{i}$, for $i=1,2$. Suppose $x$ is not adjacent to $y$. For any $\displaystyle z=\beta p_{j}\in N_{p_{j}}$ for $i\neq j$, we have $\displaystyle y+z \in U(\mathbb{Z}_{2^{m}p^{n_{1}}_{1}p^{n_{2}}_{2}})$. Now $y$ is adjacent to $y+z$, suppose not then we have following possibilities. If $\displaystyle p_{i}\mid 2y+z$, then $p_{i}\mid \beta $ and if $\displaystyle p_{j}\mid 2y+z$, then $\displaystyle p_{j}\mid \alpha$. Similarly, if $\displaystyle 2\mid 2y+z$, then $\displaystyle 2\mid \beta$, which are not possible. Since, $\displaystyle y+z,x\in U(\mathbb{Z}_{2^{m}p^{n_{1}}_{1}p^{n_{2}}_{2}})$, we have $y+z$ and $x$ is adjacent to $0$. Hence, we get $[y,y+z], [y+z,0]$ and $[0,x]$ are edges in $\displaystyle G(\mathbb{Z}_{2^{m}p^{n_{1}}_{1}p^{n_{2}}_{2}})$, which implies $d(x,y)\leq 3$.\\
(d) If $\displaystyle y\in N_{2p_{1}p_{2}}$, then $x$ is adjacent to $y$, which gives $d(x,y)=1$.\\
{\bf Case III:} If $\displaystyle x,y \in N_{U}(\mathbb{Z}_{2^{m}p^{n_{1}}_{1}p^{n_{2}}_{2}})$, then we have the following possibilities:\\
(a) If $\displaystyle x,y \in N_{2}$, then for any $\displaystyle z\in N_{p_{1}p_{2}}$,  we have $x,y$ are adjacent to $z$. Hence, we get $d(x,y)\leq 3$.\\
(b) If $\displaystyle x\in N_{2}$ and $\displaystyle y\in N_{p_{i}} $, then $x=2\alpha$ and $y=\alpha' p_{i}$. Suppose $x$ is not adjacent to $y$, then $\displaystyle p_{j}\mid x+y$, for $i\neq j$ i.e. $p_{j}\mid 2\alpha+\alpha' p_{i}$, which gives $2\alpha +\alpha' p_{i}=\beta p_{j}$. Note that $\displaystyle 2\nmid \beta$ and $\displaystyle p_{i}\nmid \beta$ this implies $\beta p_{j}\in N_{p_{j}}$ and $\displaystyle \beta p_{i}p_{j}\in N_{p_{i}p_{j}}$. First  we prove that $\displaystyle y=\alpha' p_{i}$ is adjacent to $\displaystyle z=\alpha' p_{i}+\beta p_{j}$, suppose this is not true, then we have, $\displaystyle 2\mid z+y$ i.e. $2\mid 2\alpha'p_{i}+\beta p_{j}$, which gives $2\mid \beta$ which is a contradiction. Similar, contradiction occurs, when $\displaystyle p_{i}\mid y+z$ and $\displaystyle p_{j}\mid y+z$. Now we prove that $z$ is adjacent to $u=\beta p_{i}p_{j}\in N_{p_{i}p_{j}}$, suppose it is not holds, which implies $\displaystyle 2\mid z+u$ i.e. $2\mid \alpha'p_{i}+\beta p_{j}+\beta p_{i}p_{j}$, we get $\displaystyle 2\mid \alpha' p_{i}+2\alpha +\alpha'p_{i}+\beta p_{i}$, which implies $2\mid \beta$ which is not possible. Similar contradiction occurs, when  $\displaystyle p_{i}\mid u+z$ and $\displaystyle p_{j}\mid u+z$. Finally, observe that $u=\beta p_{i}p_{j}$ is adjacent to $x=2\alpha$. Hence, we get $[y,y+z],  \ [y+z,u]$ and $[u,x]$ are edges in $\displaystyle G(\mathbb{Z}_{2^{m}p^{n_{1}}_{1}p^{n_{2}}_{2}})$, which implies that $d(x,y)\leq 3$.\\
(c) If  $\displaystyle x\in N_{2}$  and $\displaystyle y\in N_{2p_{i}}$, for $i=1,2$ then $x=2\alpha$ and $y=2\alpha' p_{i}$. First we prove that $x$ is adjacent to some, $z\in N_{p_{j}}$, for $i\neq j$. Let $\beta p_{j}\in N_{p_{2}}$ and assume that $x$ is not adjacent to $\beta p_{j}$, which implies that $\displaystyle p_{i}\mid x+\beta p_{j}$. Note that $\displaystyle \beta p_{j}\in N_{p_{2}}$, which gives $\displaystyle -\beta p_{j}\in N_{p_{j}}$ and if $x$ is not adjacent to $\displaystyle -\beta p_{j}$, then $\displaystyle p_{i}\mid x-\beta p_{j}$, this implies $\displaystyle p_{i}\mid x$, which is a contradiction. Hence, we get $\displaystyle z\in N_{p_{j}}$ such that $x$ is adjacent to $z$. Note that $\displaystyle y=2\alpha' p_{i}$ is adjacent to every $z$ in $\displaystyle N_{p_{j}}$. Hence, we get $d(x,y)\leq 3$.\\
(d) If  $\displaystyle x\in N_{2}$ and $\displaystyle y\in N_{p_{1}p_{2}}$, then $x$ is adjacent to $y$, which gives $d(x,y)=1$.\\
(e) If  $\displaystyle x\in N_{2}$ and $\displaystyle y\in N_{2p_{1}p_{2}}$, then $x=2\alpha$ and $y=2\alpha' p_{1}p_{2}$. Note that for any $\displaystyle z=\beta p_{1}p_{1}$, we have $\displaystyle x+z\in U(\mathbb{Z}_{2^{m}p^{n_{1}}_{1}p^{n_{2}}_{2}})$. Then $x$ is adjacent to $x+z$. Since $x+z$ is unit, $x+z$ is adjacent $y$. Hence, we have $[x,x+z]$ and $[x+z, y]$ are edges in $\displaystyle G(\mathbb{Z}_{2^{m}p^{n_{1}}_{1}p^{n_{2}}_{2}})$, which gives $d(x,y)\leq 3$.\\
(f) If $\displaystyle x,y\in N_{P_{i}}$, then $x=\alpha p_{i}$ and  $y=\alpha' p_{i}$ for $i=1,2$. Note that, for any $\displaystyle x\in N_{p_{i}}$ and $\displaystyle z\in N_{2p_{j}}$, for $i\neq j$, we have $x$ is adjacent to $z$, which gives $d(x,y)\leq 2$.\\
(g) If  $\displaystyle x\in N_{p_{i}}$ and $\displaystyle y\in N_{p_{j}} $, for $i\neq j$ then $x=\alpha p_{i}$ and $y=\alpha' p_{j}$. Note that $x$ is not adjacent to $y$, we prove that $x+y$ is adjacent to $x$ and $y$. Assume that $x+y$ is not adjacent to $x$, which gives $\displaystyle 2\mid 2x+y$, we get $2\mid \alpha'$ which is not possible and if $\displaystyle p_{i}\mid 2x+y$, then $p_{i}\mid \alpha'$ which is a contradiction. Similar contradiction occurs, if $\displaystyle p_{j}\mid 2x+y$. Hence we get $x$ is adjacent to $x+y$. Similarly, we can show that $y$ is adjacent to $x+y$, which gives $d(x,y)\leq 3$.\\
(h) If $\displaystyle x\in N_{p_{i}}$  and $\displaystyle y\in N_{2p_{i}} $, then $x=\alpha p_{i}$ and $y=2\alpha' p_{i}$. Let $\displaystyle \beta p_{j}\in N_{p_{j}}$, for $i\neq j$. Then $\displaystyle z=2(\beta p_{j}+\alpha p_{i} )p_{j}\in N_{2p_{j}}$, we already proved $x$ is adjacent to every $\displaystyle z\in N_{2p_{j}}$. Now we prove that $\displaystyle z$ is adjacent to $\displaystyle z+x$ and $\displaystyle z-x$. Suppose this is not hold, then $\displaystyle 2\mid 2z+x$ i.e. $2\mid\alpha p_{i}$, which gives $2\mid \alpha$ or if $\displaystyle p_{i}\mid 2z+x$ i.e. $p_{i}\mid 4(\beta p_{j}+\alpha p_{i} )p_{j}+\alpha p_{i}$, then $p_{i}\mid \beta$ or if $\displaystyle p_{j}\mid 2z+x$ i.e. $p_{j}\mid 4(\beta p_{j}+\alpha p_{i})p_{j}+\alpha p_{i}$, then  $p_{j}\mid \alpha$ which are not possible. Similarly, we can prove that $z-x$ is adjacent to $z$. Now we claim that either $z+x$ or $z-x$ is adjacent to $y$, assume that both are not true, then only possibility is $\displaystyle p_{j}\mid z+x+y$ and $\displaystyle p_{j}\mid z-x+y$, which gives $\displaystyle p_{j}\mid z+y$. Hence, we get $p_{j}\mid \alpha'$, which is a contradiction. Hence, either $z+x$ or $z-x$ is adjacent to $y$, which gives $d(x,y)\leq 3$.\\
(i) If $\displaystyle x\in N_{p_{i}}$ and $\displaystyle y\in N_{2p_{j}}$ for $i\neq j$, then $x$ is adjacent to $y$, which gives $d(x,y)=1$.\\
(j) If $\displaystyle x\in N_{p_{i}} $ and $\displaystyle y\in N_{p_{1}p_{2}}$, for $i=1,2$ then $x=\alpha p_{i}$ and $y=\alpha' p_{1}p_{2}$. We know that $y$ is adjacent to every $z=2\beta\in N_{2}$. It is enough to prove that $x$ is adjacent to some $\displaystyle 2\beta \in N_{2}$. If $\displaystyle z=2\beta \in N_{2}$ is not adjacent to $x=\alpha p_{i}$, then $p_{j}\mid z+x$, for $i\neq j$. Now note that $\displaystyle -z=-2\beta\in N_{2}$ and if $x$ is not adjacent to $-z$ then $\displaystyle p_{j}\mid x-z$ this implies $\displaystyle p_{j}\mid 2x$ which gives $p_{j}\mid \alpha$ which is a contradiction. Hence we get either $z$ or $-z$ is adjacent to $x$. Hence, $d(x,y)\leq 3$.\\
(k) If $\displaystyle x\in N_{p_{i}}$ and $\displaystyle y\in N_{2p_{1}p_{2}}$, for $i=1,2$, then using the procedure in $(e)$, we can show that $d(x,y)\leq 3$.\\
(l) If $\displaystyle x,y\in N_{2p_{i}}$, for $i=1,2$ then for any $\displaystyle z\in N_{p_{j}}$, for $i\neq j$, we have $x$ and $y$ are adjacent to $z$. Hence, $d(x,y)\leq 3$.\\
(m) If $\displaystyle x\in N_{2p_{1}}$ and $\displaystyle y\in N_{2p_{2}} $, then $x=2\alpha p_{1} $ and $y=2\alpha' p_{2}$ . Note that for any $\displaystyle \beta p_{2}\in N_{p_{2}}$, we have $\displaystyle x+\beta p_{2} $ is adjacent to $x$. Now it is enough to prove that, $y$ is adjacent to $x+\beta p_{2}$, for some $\displaystyle \beta p_{2} \in N_{p_{2}}$. Let $\displaystyle \beta p_{2}\in N_{p_{2}}$. Then $\displaystyle -\beta p_{2}\in N_{p_{2}}$, if $y$ is not adjacent to $x+\beta p_{2}$ and $x-\beta p_{2}$, then $\displaystyle p_{2}\nmid x\pm\beta p_{2}+y$ and $\displaystyle 2\nmid x\pm\beta p_{2}+y$. Hence, only possibility is that $\displaystyle p_{1}\mid x+\beta p_{2}+y$ and $\displaystyle p_{1}\mid x-\beta p_{2}+y$,  which implies $p_{1}\mid x+y$, from this we have $p_{1}\mid \alpha'$, which is a contradiction. Hence, either $x+\beta p_{2}$ or $x-\beta p_{2}$ is adjacent to $y$, which gives $d(x,y)\leq 3$.\\
(n) If $\displaystyle x\in N_{2p_{i}}$, for $i=1,2$ and $\displaystyle y\in N_{p_{1}p_{2}} $, then $x=2\alpha p_{i}$ and $y=\alpha' p_{1}p_{2}$. For any $\displaystyle 2\beta \in N_{2}$, we have $y$ is adjacent to $2\beta$, and $2\beta +y$ is adjacent to $2\beta$. Now we prove that for some $z\in N_{2}$, $x$ is adjacent to $y+z$. Let $\displaystyle 2\beta\in N_{2}$. Then  $\displaystyle -2\beta\in N_{2}$ and assume that $y+2\beta$ and $y-2\beta$ are not adjacent to $x$. Note that $\displaystyle p_{i}\nmid y+2\beta +x$ and $\displaystyle 2\nmid y+2\beta +x$, hence the only possibility is $\displaystyle p_{j}\mid y+2\beta +x$ for $i\neq j$. Similarly from above we get $\displaystyle p_{j}\mid y-2\beta +x$ this implies $p_{j}\mid y+x$, since $p_{j}\in\{p_{1},p_{2}\}$, we get $p_{j}\mid \alpha$ which is a contradiction, which gives either $y+2\beta$ is adjacent to $x$ or $y-2\beta$ is adjacent to $x$. From this, we conclude that $d(x,y)\leq 3$.\\
(o) If $\displaystyle x\in N_{2p_{i}}$ and $\displaystyle y\in N_{2p_{1}p_{2}}$, for $i=1,2$ then using the procedure in $(e)$ we can show that $d(x,y)\leq 3$.\\
(p) If $\displaystyle x,y\in N_{p_{1}p_{2}}$, then for any $\displaystyle z\in N_{2}$ we have $x$ and $y$ are adjacent to $z$. Hence we get $d(x,y)\leq 3$.\\
(q) If $\displaystyle x\in N_{p_{1}p_{2}}$ and $\displaystyle y\in N_{2p_{1}p_{2}}$, then using the procedure in $(e)$ we can show that $d(x,y)\leq 3$.\\
(r) If $\displaystyle x,y\in N_{2p_{1}p_{2}}$, then for any $\displaystyle z\in U(\mathbb{Z}_{2^{m}p^{n_{1}}_{1}p^{n_{2}}_{2}})$, we have $x$ and $y$ are adjacent to $z$, which gives $d(x,y)\leq 3$.\\
Hence, $\displaystyle G(\mathbb{Z}_{2^{m}p^{n_{1}}_{1}p^{n_{2}}_{2}})$ is connected graph and $d(x,y)\leq 3$, for all $\displaystyle x,y\in \mathbb{Z}_{2^{m}p^{n_{1}}_{1}p^{n_{2}}_{2}}$, which gives $\displaystyle \text{diam}(G(\mathbb{Z}_{2^{m}p^{n_{1}}_{1}p^{n_{2}}_{2}}))\leq 3$.
\end{proof}
\end{theorem}
\begin{corollary}\label{thm46}
Let $\displaystyle G(\mathbb{Z}_{2^{m}p^{n_{1}}_{1}p^{n_{2}}_{2}})$ be a unit graph, where $\displaystyle p_{1}$ and  $\displaystyle p_{2}$ are distinct odd primes. Then $\displaystyle g_{r}(G(\mathbb{Z}_{2^{m}p^{n_{1}}_{1}p^{n_{2}}_{2}}))=4$.
\begin{proof}
Proof follows from Theorem (\ref{thm7}) and Lemma (\ref{thm31}).
\end{proof}
\end{corollary}

\begin{theorem}\label{thm44}
Let $\displaystyle G(\mathbb{Z}_{2^{m}p^{n_{1}}_{1}p^{n_{2}}_{2}})$ be a unit graph, where $\displaystyle p_{1}$ and $\displaystyle p_{2}$  are odd primes. Then
\begin{enumerate}
  \item $\displaystyle |V|=2^{m}p^{n_{1}}_{1}p^{n_{2}}_{2}$ and $\displaystyle |E|=\frac{(2^{m}p^{n_{1}}_{1}p^{n_{2}}_{2})\phi(2^{m}p^{n_{1}}_{1}p^{n_{2}}_{2})}{2}$
  \item $\displaystyle \lambda(G(\mathbb{Z}_{2^{m}p^{n_{1}}_{1}p^{n_{2}}_{2}}))=\delta(G(\mathbb{Z}_{2^{m}p^{n_{1}}_{1}p^{n_{2}}_{2}}))=\phi(2^{m}p^{n_{1}}_{1}p^{n_{2}}_{2})$.
\end{enumerate}
\begin{proof}
\begin{enumerate}
  \item Clearly, $\displaystyle |V|=2^{m}p^{n_{1}}_{1}p^{n_{2}}_{2}$. Since, $\displaystyle 2\in N_{U}(\mathbb{Z}_{2^{m}p^{n_{1}}_{1}p^{n_{2}}_{2}})$, from Theorem (\ref{thm1}), we get $\displaystyle G(\mathbb{Z}_{2^{m}p^{n_{1}}_{1}p^{n_{2}}_{2}})$ is $\displaystyle \phi(2^{m}p^{n_{1}}_{1}p^{n_{2}}_{2})$-regular graph and hence $\displaystyle |E|=\frac{(2^{m}p^{n_{1}}_{1}p^{n_{2}}_{2})\phi(2^{m}p^{n_{1}}_{1}p^{n_{2}}_{2})}{2}$.
  \item From Lemma (\ref{thm31}) and Theorems (\ref{thm43}) and (\ref{thm4}), we get\\
   $\displaystyle \lambda(G(\mathbb{Z}_{2^{m}p^{n_{1}}_{1}p^{n_{2}}_{2}}))=\delta(G(\mathbb{Z}_{2^{m}p^{n_{1}}_{1}p^{n_{2}}_{2}}))=\phi(2^{m}p^{n_{1}}_{1}p^{n_{2}}_{2})$.
\end{enumerate}
\end{proof}
\end{theorem}

\begin{theorem} Let $\displaystyle G(\mathbb{Z}_{n})$ be a unit graph, where $n=p^{n_{1}}_{1}p^{n_{2}}_{2}p^{n_{3}}_{3}$ and $p_{1}, \ p_{2}  \ \& \ p_{3}$ are distinct primes. Let $H$ be a $|V|\times |E|$ incidence matrix of $G(\mathbb{Z}_{n})$.
    \begin{enumerate}
      \item If $\displaystyle 2\in U(\mathbb{Z}_{n})$, then binary code generated by $H$ is a $\displaystyle C_{2}(H)=\left[\frac{(n-1)\phi(n)}{2},n-1,\phi(n)-1\right]_{2}$ code over finite field $\mathbb{F}_{2}$.
      \item If $\displaystyle 2\in N_{U}(\mathbb{Z}_{n})$, then for any odd prime $q$, the $q$-ary code generated by $H$ is a $\displaystyle C_{q}(H)=\left[\frac{n\phi(n)}{2},n-1,\phi(n)\right]_{q}$ code over finite field $\mathbb{F}_{q}$.
    \end{enumerate}
\begin{proof}
\begin{enumerate}
  \item Let $\displaystyle 2\in U(\mathbb{Z}_{n})$ then $p_{1},p_{2}$ and $p_{3}$ are odd primes. By Theorem (\ref{thm41}), $\displaystyle G(\mathbb{Z}_{n})=(V,E)$ is a connected graph and hence by Theorem (\ref{thm2}) binary code generated by $H$ is
      $\displaystyle C_{2}(H)=[|E|,|V|-1,\lambda(G(\mathbb{Z}_{p^{n_{1}}_{1}p^{n_{2}}_{2}p^{n_{3}}_{3}}))]_{2}$. Now from Theorem (\ref{thm42}), we get $\displaystyle |E|= \frac{(n-1)\phi(n)}{2}, \ |V|-1=n-1$ and $\displaystyle \lambda(G(\mathbb{Z}_{n}))=\phi(n)-1$.
  \item Let $\displaystyle 2\in N_{U}(\mathbb{Z}_{n})$. Then we have $p_{i}=2$, for some $i$. By Theorem (\ref{thm44}) and Lemma (\ref{thm31}), $\displaystyle G(\mathbb{Z}_{n})=(V,E)$ is a connected bipartite graph and hence by Theorems (\ref{thm2}) and (\ref{thm5}), $q$-ary code generated by $H$ is $\displaystyle C_{q}(H)=[|E|,|V|-1,\lambda(G(\mathbb{Z}_{n}))]_{q}$. Now from Theorem (\ref{thm44}), the result follows.
\end{enumerate}
\end{proof}
\end{theorem}

\begin{corollary}
Let $C_{q}(H)$ and $C_{2}(H)$ denote the linear codes generated from incidence matrices of $\displaystyle G(\mathbb{Z}_{2^{m}p^{n_{1}}_{1}p^{n_{2}}_{2}})$ and $\displaystyle G(\mathbb{Z}_{p^{n_{1}}_{1}p^{n_{2}}_{2}p^{n_{3}}_{3}})$. Then
\begin{enumerate}
  \item Dual of code $C_{2}$ is $\displaystyle C^{\perp}_{2}=\left[\frac{(n-1)\phi(n)}{2},\frac{(n-1)[\phi(n)-2]}{2},3\right]_{2}$,\\
   where $\displaystyle n=p^{n_{1}}_{1}p^{n_{2}}_{2}p^{n_{3}}_{3}$.
  \item Dual of code $C_{q}$ is $\displaystyle C^{\perp}_{q}=\left[\frac{n\phi(n)}{2},\frac{n(\phi(n)-2)+2}{2},4\right]_{q}$,\\
   where $\displaystyle n=2^{m}p^{n_{1}}_{1}p^{n_{2}}_{2}$.
\end{enumerate}
\begin{proof}
\begin{enumerate}
\item From Theorem (\ref{thm8}), $\text{diam}(\displaystyle C^{\perp}_{2})=\frac{(n-1)[\phi(n)-2]}{2}$. By Theorem (\ref{thm6}) and Corollary  (\ref{thm45}), $\displaystyle d(C^{\perp}_{2})=3$.
  \item Proof follows from Theorems (\ref{thm8}) \& (\ref{thm6}) and Corollary (\ref{thm46}).
\end{enumerate}
\end{proof}
\end{corollary}
\begin{theorem}
Let $G(\mathbb{Z}_{n})$ be a unit graph, where $n$ is any positive integer.
\begin{enumerate}
  \item If $2\in U(\mathbb{Z}_{n})$, then $\displaystyle |E|=\frac{(n-1)\phi(n)}{2}$.
  \item If $2\in N_{U}(\mathbb{Z}_{n})$, then $\displaystyle |E|=\frac{n\phi(n)}{2}$.
\end{enumerate}
\begin{proof}
\begin{enumerate}
  \item Let $\displaystyle 2\in U(\mathbb{Z}_{n})$. By Theorem (\ref{thm1}), $\displaystyle \text{deg}(x)=\phi(n)-1$ for all $\displaystyle x\in U(\mathbb{Z}_{n})$ and $\displaystyle \text{deg}(x)=\phi(n)$ for all $\displaystyle x\in N_{U}(\mathbb{Z}_{n})$. Thus,
			\begin{align*}
				|E| &= \frac{\sum_{x\in V}\text{deg}(x)}{2}\\
				&= \frac{\sum_{x\in U(\mathbb{Z}_{n})}\text{deg}(x)+\sum_{x\in N_{U}(\mathbb{Z}_{n})}\text{deg}(x)}{2}\\
				&= \frac{\phi(n)(\phi(n)-1)+ (n-\phi(n))\phi(n)}{2}\\
				|E| &= \frac{(n-1)\phi(n)}{2}.
			\end{align*}
\item Let $\displaystyle 2\in N_{U}(\mathbb{Z}_{n})$. Then from Theorem (\ref{thm1}), we get $\displaystyle G(\mathbb{Z}_{n})$ is $\displaystyle \phi(n)$-regular graph.  Hence, $\displaystyle |E|=\frac{n\phi(n)}{2}$.
\end{enumerate}
\end{proof}
\end{theorem}
Bases on our results for unit graphs, in Theorems (\ref{thm22}), (\ref{thm26}), (\ref{thm32}), (\ref{thm34}), (\ref{thm41}) and (\ref{thm43}), we state following conjectures for any natural number $n$: \\
{\bf Conjecture I:}
Let $G(\mathbb{Z}_{n})$ be a unit graph. Then $G(\mathbb{Z}_{n})$ is connected graph and
\begin{enumerate}
  \item If $2\in U(\mathbb{Z}_{n})$, then $\text{diam}(G(\mathbb{Z}_{n}))\leq 2$.
  \item If $2\in N_{U}(\mathbb{Z}_{n})$, then $\text{diam}(G(\mathbb{Z}_{n}))\leq 3$.
\end{enumerate}
{\bf Conjecture II:}
Let $\displaystyle G(\mathbb{Z}_{n})$ be a unit graph and $H$ be a $|V|\times |E|$ incidence matrix of $G(\mathbb{Z}_{n})$.
    \begin{enumerate}
      \item If $\displaystyle 2\in U(\mathbb{Z}_{n})$, then binary code generated by $H$ is a $\displaystyle C_{2}(H)=\left[\frac{(n-1)\phi(n)}{2},n-1,\phi(n)-1\right]_{2}$ code over finite field $\mathbb{F}_{2}$.
      \item If $\displaystyle 2\in N_{U}(\mathbb{Z}_{n})$, then for any odd prime $q$, the $q$-ary code generated by $H$ is a $\displaystyle C_{q}(H)=\left[\frac{n\phi(n)}{2},n-1,\phi(n)\right]_{q}$ code over finite field $\mathbb{F}_{q}$.
    \end{enumerate}
\section{Conclusion}
In this paper, we generate $q$-ary linear codes, for any prime $q$, from incidence matrices of unit graphs $\displaystyle G(\mathbb{Z}_{n})$. Moreover, we found parameters and dual of constructed codes over finite filed $\displaystyle \mathbb{F}_{q}$. In this article, we consider $n$ as product of power of three distinct primes. Examine the permutation decoding techniques, covering radius of constructed codes and one can construct linear codes from unit graph over different commutative rings is the further scope to work.

\end{document}